\documentclass[aap,citesort,MSNbibl,dvips]{arximspdf}
\usepackage{mathbh}
\usepackage{graphicx}

\doi{10.1214/10-AAP682}
\volume{20}
\issue{6}
\pubyear{2010}
\firstpage{2178}
\lastpage{2203}

\makeatletter

\newcommand{\ball}{\overline{B}}
\newcommand{\E}{\mathbb{E}}
\renewcommand{\P}{\mathbb{P}}

\newcommand{\defeq}{:=}
\newcommand{\eqdef}{=:}
\newcommand{\ud}{{d}}
\newcommand{\R}{\mathbb{R}}
\newcommand{\F}{\mathcal{F}}
\newcommand{\trace}{\operatorname{trace}}
\newcommand{\uarg}{ \cdot}

\newproclaim{definition}{Definition}
\newtheorem{corollary}[definition]{Corollary}
\newtheorem{theorem}[definition]{Theorem}
\newtheorem{lemma}[definition]{Lemma}
\newtheorem{proposition}[definition]{Proposition}
\newproclaim{remark}[definition]{Remark}

\makeatother

\begin{document}
\begin{frontmatter}

\title{On the ergodicity of the adaptive Metropolis algorithm
on unbounded domains}
\runtitle{On the ergodicity of adaptive Metropolis}

\begin{aug}
\author[A]{\fnms{Eero} \snm{Saksman}\thanksref{t1}\ead[label=e1]{eero.saksman@helsinki.fi}} and
\author[B]{\fnms{Matti} \snm{Vihola}\corref{}\thanksref{t2}\ead[label=e2]{matti.vihola@iki.fi}
\ead[label=u2,url]{http://iki.fi/mvihola}}
\runauthor{E. Saksman and M. Vihola}
\affiliation{University of Helsinki and University of Jyv\"{a}skyl\"{a}}
\address[A]{Department of Mathematics and Statistics\\
P.O. Box 68\\
FI-00014 University of Helsinki\\
Finland \\
\printead{e1}}
\address[B]{Department of Mathematics and Statistics\\
P.O. Box 35 (MaD)\\
FI-40014 University of Jyv\"{a}skyl\"{a}\\
Finland \\
\printead{e2}\\
\printead{u2}}
\end{aug}

\thankstext{t1}{Supported by the Academy of Finland, Projects 113826
and 118765,
and by the Finnish Centre of Excellence in Analysis and Dynamics
Research.}

\thankstext{t2}{Supported by the Academy of Finland, Projects 110599
and 201392,
by the Finnish Academy of Science and Letters, Vilho, Yrj\"{o} and
Kalle V\"{a}is\"{a}l\"{a} Foundation, by the Finnish Centre of
Excellence in Analysis and Dynamics Research and by the Finnish
Graduate School in Stochastics and Statistics.}

\received{\smonth{6} \syear{2008}}
\revised{\smonth{2} \syear{2009}}

%
\begin{abstract}
This paper describes sufficient conditions to ensure the correct
ergodicity of the Adaptive Metropolis (AM) algorithm of Haario, Saksman
and Tamminen [\textit{Bernoulli} \textbf{7} (2001) 223--242] for target
distributions with a noncompact support. The conditions ensuring a
strong law of large numbers require that the tails of the target
density decay super-exponentially and have regular contours. The result
is based on the ergodicity of an auxiliary process that is sequentially
constrained to feasible adaptation sets, independent estimates of the
growth rate of the AM chain and the corresponding geometric drift
constants. The ergodicity result of the constrained process is obtained
through a modification of the approach due to Andrieu and Moulines
[\textit{Ann. Appl. Probab.} \textbf{16} (2006) 1462--1505].
\end{abstract}

%

\end{frontmatter}

\section{Introduction} 

The Markov chain Monte Carlo (MCMC) method, first proposed by
\cite{metropolis}, is a commonly used device for numerical
approximation of
integrals of the type
\[
\pi(f) = \int f(x) \pi(x) \,\ud x,
\]
where $\pi$ is a probability density function. Intuitively, the method is
based on producing a sample $(X_k)_{k=1}^n$ of random variables from the
distribution $\pi$ defines. The integral $\pi(f)$ is
approximated with the average $I_n\defeq n^{-1}\sum_{k=1}^n f(X_k)$. In
particular,
the random variables $(X_k)_{k=1}^n$ are a realization of a Markov chain,
constructed so that the chain has $\pi$ as the unique invariant
distribution.

One of the most commonly applied constructions of such a chain in
$\R^d$ is to let $X_0\equiv x_0$ with some fixed point $x_0\in\R^d$,
and recursively for $n\ge1$:
\begin{enumerate}
\item simulate $Y_{n} = X_{n-1} + U_n$, where $U_n$ is an independent random
variable distributed according to some symmetric proposal distribution $q$,
for example, a~zero-mean Gaussian, and
\item with probability $\min\{1,
\pi(Y_{n})/\pi(X_{n-1})\}$, the proposal is accepted and
$X_{n} = Y_{n}$; otherwise the proposal is rejected and $X_{n} = X_{n-1}$.
\end{enumerate}
This symmetric random-walk Metropolis algorithm is often efficient enough,
even in a relatively complex and high-dimensional situation, provided that
the proposal distribution $q$ is selected properly. Finding a good proposal
for a particular problem can, however, be a difficult task.

Recently, there has been a number of publications describing different
adaptation techniques aiming to find a good proposal automatically
\cite
{saksman-am,andrieu-robert,atchade-rosenthal,andrieu-moulines,roberts-rosenthal}
(see also the review article \cite{andrieu-thoms}).
It has been a common practice to perform trial runs, and determine the
proposal from the outcome. The recently proposed methods are different in
that they adapt on-the-fly, continuously during the estimation run. In
this paper, we focus on the forerunner of these methods, the Adaptive
Metropolis (AM) algorithm \cite{saksman-am}, which is a random-walk
Metropolis sampler with a Gaussian proposal $q_v$ having a covariance
$v$. The proposal covariance $v$ is updated continuously during the run,
according to the history of the chain. In general, such an adaptation
may, if carelessly implemented, destroy the correct ergodicity properties,
that is, that $I_n$ does not converge to $\pi(f)$ as $n\to\infty$
(see, e.g.,
\cite{roberts-rosenthal}). For practical considerations
of the AM algorithm, the reader may consult
\cite{saksman-jrstatsoc,roberts-rosenthal-examples}.

In the original paper \cite{saksman-am} presenting the AM algorithm,
the first ergodicity result for such adaptive algorithms was obtained. More
precisely, a strong law of large numbers was proved for bounded
functionals, when the
algorithm is run on a compact subset of $\R^d$. After that, several authors
have obtained more general conditions under which an adaptive MCMC process
preserves the correct ergodicity properties. Andrieu and Robert
\cite{andrieu-robert} established the connection between adaptive MCMC and
stochastic approximation, and proposed a general framework for adaptation.
Atchad\'{e} and Rosenthal \cite{atchade-rosenthal} developed further the
technique of \cite{saksman-am}. Andrieu and Moulines
\cite{andrieu-moulines} made important progress by generalizing
the Poisson equation and martingale approximation techniques
to the adaptive setting. They proved the ergodicity and a central limit
theorem for a class of adaptive MCMC schemes. Roberts and Rosenthal
\cite{roberts-rosenthal}
use an interesting approach based on coupling to show a weak law of large
numbers. However, in
the case of AM, all the techniques essentially assume that the adapted
parameter is constrained to a predefined compact set, or do not present
concrete verifiable conditions. The only result to
overcome this assumption is the one by Andrieu and Moulines
\cite{andrieu-moulines}. Their
result, however, requires a modification of the algorithm, including
additional re-projections back to some fixed compact set.

This paper describes sufficient conditions under which the AM algorithm
preserves the correct ergodicity properties, and $I_n\to\pi(f)$ almost
surely as $n\to\infty$ for any function $f$ that is bounded on
compact sets
and grows at most exponentially as $\|x\|\to\infty$. Our main result
(Theorem \ref{thm:am-ergodicity}) holds
for the original AM process (without re-projections) having a target
distribution supported on~$\R^d$. Essentially, the target density $\pi$
must have asymptotically lighter tails than $\pi(x) = c e^{-\|x\|^p}$ for
some $p>1$, and for large enough $\|x\|$, the sets $A_x=\{y\in
\R^d\dvtx\pi(y)\ge\pi(x)\}$ must have uniformly regular contours. Our
assumptions are very close to the well-known conditions proposed by Jarner
and Hansen \cite{jarner-hansen} to ensure the geometric convergence of a
(nonadaptive) Metropolis process.
By the techniques of this paper, one may also establish a central limit
theorem (see Theorem \ref{thm:am-clt}).

The ergodicity results for the AM process rely on three main contributions.
First, in Section \ref{sec:notations}, we describe an adaptive MCMC
framework, in
which the adaptation parameter is constrained at each time to a feasible
adaptation set. In Section \ref{sec:main-andrieu}, we prove
a strong law of large numbers for such a process,
through a modification of the technique of
Andrieu and Moulines \cite{andrieu-moulines}. Second, we propose
an independent estimate for the growth rate of a process satisfying a
general drift condition in Section \ref{sec:grow-bound}. Third, in
Section \ref{sec:am-ergodicity}, we provide an estimate for constants of
geometric drift for a symmetric random-walk Metropolis process, when the
target distribution has super-exponentially decaying tails with
regular contours.

The paper is essentially self-contained, and assumes little background
knowledge. Only the basic martingale theory is needed to
follow the argument, with the exception of Theorem
\ref{th:m-t-computable} by Meyn and Tweedie
\cite{meyn-tweedie-computable}, restated in Appendix
\ref{sec:proof-lemma-drift-to-conv}. Even though we consider only the AM
algorithm, our techniques apply also to many other adaptive MCMC
schemes of similar type.


\section{General framework and notation}
\label{sec:notations} 

We consider an adaptive Markov chain Monte Carlo (MCMC) chain evolving in
space $\mathbb{X}\times\mathbb{S}$, where $\mathbb{X}$ is the state
space of the ``MCMC'' chain $(X_n)_{n\ge0}$ and
the adaptation parameter $(S_n)_{n\ge0}$ evolves in
$\mathbb{S}\subset\overline{\mathbb{S}}$, where $\overline{\mathbb
{S}}$ is a
separable normed
vector space.
We assume an underlying probability space $(\Omega,\F_\Omega,\P)$,
and denote the expectation with respect to $\P$ by $\E$.
The natural filtration of the chain is denoted with $\F\defeq(\F
_k)_{k\ge
0}\subset\F_{\Omega}$ where $\F_k\defeq\sigma(X_j,S_j\dvtx0\le
j\le k)$.
We also assume that we are given an increasing sequence $K_0\subset
K_1\subset\cdots\subset K_n \subset\mathbb{S}$ of subsets of the
adaptation parameter space $\mathbb{S}$.
The random variables
$(X_n,S_n)_{n\ge0}$ form a stochastic chain, starting from
$S_0\equiv s_0\in K_0\subset\mathbb{S}$ and $X_0 \equiv x_0\in
\mathbb{X}$,
and for
$n\ge0$, satisfying the following
recursion:
%
%
\begin{eqnarray}
\label{eq:state-eq}
X_{n+1} &\sim& P_{S_n}(X_n, \uarg), \\
\label{eq:adapt-eq}
S_{n+1} &=& \sigma_{n+1} (S_n, \eta_{n+1}H(S_n, X_{n+1}) ),
\end{eqnarray}
where $P_s$ is a transition probability for each $s\in\mathbb{S}$,
$H\dvtx\mathbb{S}\times\mathbb{X}\to\overline{\mathbb{S}}$ is an
adaptation
function, and
$(\eta_n)_{n\ge1}$ is a decreasing sequence of adaptation step sizes
$\eta_n\in(0,1)$.
The functions $\sigma_n\dvtx\mathbb{S}\times\overline{\mathbb
{S}}\to
\mathbb{S}$
are defined as
\[
\sigma_n(s,s') \defeq\cases{
s+s',&\quad if $s+s'\in K_n$, \cr
s,&\quad otherwise.}
\]
Thus, $\sigma_n$ ensures that $S_n$ lies in $K_n$ for each $n\ge0$.
The recursion (\ref{eq:adapt-eq}) can also be considered as
constrained Robbins--Monro stochastic approximation (see
\cite{andrieu-moulines,sa-verifiable} and references therein).

Let $V\dvtx\mathbb{X}\to[1,\infty)$ be a function. We define a
$V$-norm of a
function $f$ as
\[
\|f\|_V \defeq\sup_{x\in\mathbb{X}} \frac{|f(x)|}{V(x)}.
\]
As usual, we denote the integration of a function $f$ with respect to a (signed)
measure $\mu$ as $\mu(f)\defeq\int f( x) \mu(\ud x)$, and define $Pf(x)
\defeq\int f(y) P(x,\ud y)$ for a transition probability $P$.
The $V$-norm of a signed measure is defined as
\[
\|\mu\|_V \defeq{\sup_{|f|\le V}} |\mu(f)|.
\]
The indicator function of a set $A$ is denoted as $\mathbh{1}_{A}(x)$ and
equals one if $x\in A$ and zero otherwise. In addition, we use the notation
$a\vee b\defeq\max\{a,b\}$ and $a \wedge b\defeq\min\{a,b\}$.

Finally, we define the following regularity property for a family of
functions $\{f_s\}_{s\in\mathbb{S}}$.
%
\begin{definition} \label{def:poly-lip} 
Suppose $V\dvtx\mathbb{X}\to[1,\infty)$.
Given an increasing sequence of subsets $K_n\subset\mathbb{S}$, $n\ge1$,
we say that
a family of functions $\{f_s\}_{s\in\mathbb{S}}$, with $f_s\dvtx
\mathbb{X}\to\R
$, is
$(K_n,V)$-polynomially Lipschitz with constants $c\ge
1,\varepsilon\ge0$, if for all $s,s'\in K_n$, we have
\[
\|f_s\|_V \le c n^\varepsilon\quad \mbox{and}\quad
\|f_s - f_{s'}\|_V \le c
n^\varepsilon|s-s'|.
\]
\end{definition}
%


\section{Ergodicity of sequentially constrained adaptive MCMC}
\label{sec:main-andrieu} 

This section contains general ergodicity results for a sequentially
constrained process defined in Section \ref{sec:notations}. These results
can be seen auxiliary to our results on Adaptive Metropolis in Section
\ref{sec:am-ergodicity}, but may be applied to other
adaptive MCMC methods as well.

Suppose that the adaptation algorithm has the form given in
(\ref{eq:state-eq}) and (\ref{eq:adapt-eq}), and the following
assumptions are satisfied for some $c\ge1$ and $\varepsilon\ge0$.
\begin{enumerate}[(A1)]
\item[(A1)]\hypertarget{a:invariance} 
For each $s\in\mathbb{S}$, the transition probability $P_s$ has $\pi
$ as
the unique invariant
distribution. 
\item[(A2)]\hypertarget{a:uniform-drift-mino} 
For each $n\ge1$, the following uniform drift and minorization
condition holds for all $s\in K_n$:
%
%
\begin{eqnarray}
\label{eq:drift-ineq}
P_s V(x) &\le&\lambda_n V(x) + b_n \mathbh{1}_{C_n}(x)
\qquad \forall x\in\mathbb{X},\\
\label{eq:mino-ineq}
P_s(x,A) &\ge&\delta_n \nu_s(A)\qquad \forall x\in C_n, \forall
{A}\subset\mathbb{X},
\end{eqnarray}
where $C_n\subset\mathbb{X}$ is a subset (a minorization set),
$V\dvtx\mathbb{X}\to[1,\infty)$ is a drift function such that $\sup
_{x\in C_n}V(x)
\le b_n$, and $\nu_s$ is a probability measure on $\mathbb{X}$,
concentrated on $C_n$.
Furthermore, the
constants $\lambda_n\in(0,1)$ and $b_n\in(0,\infty)$
are increasing, and $\delta_n\in(0,1]$ is decreasing with respect to $n$,
and they are polynomially bounded so that
\[
(1-\lambda_n)^{-1}
\vee\delta_n^{-1}
\vee b_n
\le c n^{\varepsilon}.
\]
%
%
\item[(A3)]\hypertarget{a:kernel-lip} 
For all $n\ge1$ and any $r\in(0,1]$, there is $c'=c'(r)\ge1$ such that
for all $s,s'\in K_n$,
\[
\|P_s f - P_{s'} f \|_{V^r}
\le c' n^{\varepsilon} \|f\|_{V^r} | s - s'| .
\]
%
%
\item[(A4)]\hypertarget{a:adapt-bound} 
There is a $\beta\in[0,1/2]$
such that for all $n\ge1$, $s\in K_n$ and $x\in\mathbb{X}$
\[
|H(s,x)| \le c n^\varepsilon V^\beta(x).
\]
%
\end{enumerate}
%
%
\begin{theorem}\label{thm:slln-restricted} 
Assume \textup{\hyperlink{a:invariance}{(A1)}--\hyperlink{a:adapt-bound}{(A4)}} hold
and let $f$ be a function
with $\|f\|_{V^\alpha} < \infty$ for some $\alpha\in(0,1-\beta)$.
Assume $\varepsilon< \kappa_*^{-1}
[(1/2) \wedge(1-\alpha-\beta) ]$, where $\kappa_*\ge1$ is an
independent constant, and that
$\sum_{k=1}^\infty k^{\kappa_*\varepsilon-1} \eta_k < \infty$.
Then
%
%
\begin{equation} \label{eq:slln}
\frac{1}{n} \sum_{k=1}^n f(X_k)
\stackrel{n\to\infty}{\hbox to 1cm{\rightarrowfill}}
\pi(f)\qquad
\mbox{almost surely.}
\end{equation}
\end{theorem}
%

The proof of Theorem \ref{thm:slln-restricted} is postponed to the end of
this section. We start by the following lemma, whose proof is given in
Appendix \ref{sec:proof-lemma-drift-to-conv}. It shows that if we have
polynomially worse bounds for drift and minorization constants, then the
speed of geometric convergence can get only polynomially worse.
%
\begin{lemma} \label{lemma:drift-to-conv} 
Suppose \textup{\hyperlink{a:uniform-drift-mino}{(A2)}} holds.
Then, one has for $r\in(0,1]$ that for all $s\in K_n$ and $k\ge1$,
\[
\|P_s^k(x,\uarg)-\pi(\uarg)\|_{V^r} \le V^r(x) L_n \rho_n^k
\]
with bound
\[
L_n \vee(1-\rho_n)^{-1} \le c_2 n^{\kappa_2\varepsilon},
\]
where $\kappa_2>0$ is an independent constant,
and $c_2=c_2(c,r)\ge1$.
\end{lemma}

Observe that the statement in Lemma \ref{lemma:drift-to-conv} entails that
any function $\|f\|_{V} <\infty$ is integrable with respect to
the measures $\pi$ and $P_s^k(x,\uarg)$, for all $x\in\mathbb{X}$,
$k\ge1$
and $s\in\bigcup_{n\ge0} K_n$.
The next three results are modified from Proposition 3, Lemma 5 and
Proposition 6 of \cite{andrieu-moulines}, respectively. The first one
bounds the regularity of the solutions $\hat{f}_s$ of the Poisson equation
%
%
\begin{equation}\label{eq:poisson}
\hat{f}_s - P_s \hat{f}_s = f_s - \pi(f_s)
\end{equation}
for a polynomially Lipschitz family of functions.
%
\begin{proposition}
\label{prop:poisson-regularity} 
Suppose that \textup{\hyperlink{a:invariance}{(A1)}--\hyperlink{a:kernel-lip}{(A3)}} hold, and the
family of
functions $\{f_s\}_{s\in\mathbb{S}}$ is $(K_n,V^r)$-polynomially Lipschitz
with constants $(c,\varepsilon)$, for some $r\in(0,1]$.
There is an independent constant $\kappa_3>0$ and a constant
\mbox{$c_3 = c_3(c,c',r)\ge1$}, such that:

\begin{longlist}
\item
The family $\{P_s f_s\}_{s\in\mathbb{S}}$ is
$(K_n, V^r)$-polynomially Lipschitz with constants
$(c_3,\kappa_3\varepsilon)$.

\item
Define, for any
$s\in\mathbb{S}$, the function
%
%
\begin{equation}\label{eq:poisson-solution}
\hat{f}_s \defeq\sum_{k=0}^\infty[P_s^k f_s - \pi(f_s) ].
\end{equation}
Then, $\hat{f}_s$ solves the Poisson equation (\ref{eq:poisson}), and
the families $\{\hat{f}_s\}_{s\in\mathbb{S}}$ and $\{P_s
\hat{f}_s\}_{s\in\mathbb{S}}$ are $(K_n, V^r)$-polynomially Lipschitz
with constants $(c_3, \kappa_3\varepsilon)$. In other words,
%
%
\begin{eqnarray}
\label{eq:poisson-reg1}
\|\hat{f}_s\|_{V^r} + \|P_s\hat{f}_s\|_{V^r}
&\le& c_3 n^{\kappa_3\varepsilon},
\\
\label{eq:poisson-reg2}
\|\hat{f}_s- \hat{f}_{s'}\|_{V^r}
+ \|P_s\hat{f}_s - P_{s'}\hat{f}_{s'}\|_{V^r}
&\le& c_3 n^{\kappa_3\varepsilon}
|s-s'|
\end{eqnarray}
for all $s,s'\in K_n$.
\end{longlist}
\end{proposition}
%
%
\begin{pf} 
Throughout the proof, suppose $s,s'\in K_n$.

The part (i) follows easily from Lemma
\ref{lemma:drift-to-conv}, since
\begin{eqnarray*}
\|P_s f_s\|_{V^r} &\le&
\|P_s f_s - \pi(f_s)\|_{V^r}
+ |\pi(f_s)|\le[c_2 n^{\kappa_2 \varepsilon} + \pi(V^r)] \|
f_s\|_{V^r}, \\
\|P_s f_s - P_{s'} f_{s'}\|_{V^r}
&\le&
\|(P_s - P_{s'})f_s\|_{V^r} + \|P_{s'}(f_s-f_{s'})\|_{V^r}
\\
&\le& c' n^{\varepsilon} \|f_s\|_{V^r} |s-s'|
+ \tilde{c} n^{\kappa_2\varepsilon}\|f_s-f_{s'}\|_{V^r}
\le\tilde{c} n^{(\kappa_2+1)\varepsilon} |s-s'|.
\end{eqnarray*}

Consider then (ii). Estimate (\ref{eq:poisson-reg1}) follows
by the definition of $\hat{f}_s$ and
Lemma~\ref{lemma:drift-to-conv},
\begin{eqnarray*}
\|\hat{f}_s\|_{V^r}
&\le&\sum_{k=0}^\infty\|P_s^k f_s - \pi(f_s)\|_{V^r}
\le L_n \|f_s\|_{V^r} \sum_{k=0}^\infty\rho_n^k\\
&=&\frac{L_n}{1-\rho_n} \|f_s\|_{V^r}
\le(c_2 n^{\kappa_2\varepsilon})^2 c n^{\varepsilon}
= c_2^2 c n^{(2\kappa_2 + 1)\varepsilon}.
\end{eqnarray*}
The above bound clearly applies also to $\|P_s \hat{f}_s\|_{V^r}$,
and the convergence implies that $\hat{f}_s$ solves (\ref{eq:poisson}).

For (\ref{eq:poisson-reg2}), define an auxiliary transition probability
by setting
$\Pi(x,{A})\defeq\pi(A)$, and write
\[
P_s^k f - P_{s'}^k f
= \sum_{j=0}^{k-1}
(P_s^j - \Pi)(P_s-P_{s'})[P_{s'}^{k-j-1} f - \pi(f) ]
\]
since $\pi P_s = \pi$ for all $s$. By Lemma \ref{lemma:drift-to-conv} and
assumption \hyperlink{a:kernel-lip}{(A3)}, we have for all $s,s'\in K_n$ and
$j\ge0$
\begin{eqnarray*}
&&\| (P_s^j - \Pi)(P_s-P_{s'})
[P_{s'}^{k-j-1} f - \pi(f) ] \|_{V^r} \\
&&\qquad\le L_n \rho_n^j
\| (P_s-P_{s'})
[P_{s'}^{k-j-1} f - \pi(f) ]\|_{V^r} \\
&&\qquad\le L_n \rho_n^j
c' n^{\varepsilon} |s-s'|
\|P_{s'}^{k-j-1} f - \pi(f)\|_{V^r} \\
&&\qquad\le L_n^2 \rho_n^{k-1},
\end{eqnarray*}
which gives that
%
%
\begin{equation} \label{eq:kernel-diff-estimate}
\| P_s^k f - P_{s'}^k f\|_{V^r} \le k
L_n^2 \rho_n^{k-1}
c' n^{\varepsilon} |s-s'| \|f\|_{V^r}.
\end{equation}
Write then
\[
\hat{f}_s - \hat{f}_{s'}
= \sum_{k=0}^\infty[ P_s^k f_s - P_{s'}^k f_s ]
- \sum_{k=0}^\infty[ P_{s'}^k(f_{s'}-f_s)-\pi(f_{s'}-f_s) ].
\]
By Lemma \ref{lemma:drift-to-conv} and estimate
(\ref{eq:kernel-diff-estimate}) we have
\begin{eqnarray*}
\|\hat{f}_s - \hat{f}_{s'} \|_{V^r}
&\le& L_n^2
c' n^{\varepsilon} |s-s'|
\Biggl(\sum_{k=0}^\infty k\rho_n^{k-1} \Biggr) \|f_s\|_{V^r}
+ L_n \Biggl(\sum_{k=0}^\infty\rho_n^k \Biggr) \|f_{s'}-f_s\|_{V^r}
\\
&\le&[
L_n^2 c' n^{\varepsilon} (1-\rho_n)^{-2} c n^{\varepsilon}
+ L_n (1-\rho_n)^{-1} c n^{\varepsilon}
] |s-s'|
\\
&\le&[
(c_2 n^{\kappa_2\varepsilon})^2
c' n^{\varepsilon} (c_2 n^{\kappa_2\varepsilon})^{2} c
n^{\varepsilon}
+ (c_2 n^{\kappa_2\varepsilon}) (c_2 n^{\kappa_2\varepsilon}) c
n^{\varepsilon}
] |s-s'|
\\
&\le&2 c_2^4 c' c n^{(4\kappa_2+2)\varepsilon} |s-s'|.
\end{eqnarray*}
The same bound applies, with a similar argument, to
$P_s \hat{f}_s - P_{s'}\hat{f}_{s'}$.
\end{pf}
%
%
%
\begin{lemma}
\label{lemma:stability} 
Assume that \textup{\hyperlink{a:uniform-drift-mino}{(A2)}} holds.
Then, for all $r\in[0,1]$, any sequence $(a_n)_{n\ge1}$ of
positive numbers, and $(x_0,s_0)\in\mathbb{X}\times K_0$, we have that
%
%
\begin{eqnarray}\label{eq:stability1}
\E[V^r(X_k) ]
&\le& c_4^r k^{2r\varepsilon} V^r(x_0),
\\
\label{eq:stability2}
\E\Bigl[\max_{m \le j\le k} (a_j V(X_j))^r \Bigr]
&\le& c_4^r \Biggl(
\sum_{j=m}^k a_j j^{2\varepsilon}
\Biggr)^r V^r(x_0),
\end{eqnarray}
where the constant $c_4$ depends only on $c$.
\end{lemma}
%
%
\begin{pf} 
For $(x_0,s_0)\in\mathbb{X}\times K_0$ and $k\ge1$, we can apply
the drift inequality (\ref{eq:drift-ineq}) and the monotonicity of
$\lambda_k$ and $b_k$ to obtain
%
%
\begin{eqnarray} \label{eq:drift-estimate}\quad
\mathbb{E} [V(X_k) ]
&=& \mathbb{E} [\mathbb{E}[V(X_k) | \F_{k-1}] ] =
\mathbb{E} [P_{S_{k-1}}V(X_{k-1}) ]\nonumber\\
&\le&\lambda_k \mathbb{E} [V(X_{k-1}) ] + b_k
\le\cdots\le\lambda_k^k V(x_0) + b_k\sum_{j=0}^{k-1} \lambda
_k^j\\
&\le&\Biggl(1+b_k\sum_{j=0}^\infty\lambda_k^j\Biggr)V(x_0)
\le(1+c^2k^{2\varepsilon})V(x_0)
\le c_4 k^{2\varepsilon} V(x_0).\nonumber
\end{eqnarray}
This estimate
with Jensen's inequality yields for $r\in[0,1]$ that
\[
\mathbb{E} [V^r(X_k) ]
\le(\mathbb{E} [V(X_k) ] )^r
\le c_4^r k^{2r\varepsilon} V^r(x_0).
\]
Similarly, we have
\begin{eqnarray*}
\mathbb{E} \Bigl[\max_{m\le j\le k} (a_j V(X_j))^r \Bigr]
&\le&\Bigl(
\mathbb{E} \Bigl[\max_{m\le j\le k} a_j V(X_j) \Bigr] \Bigr)^r\\
&\le&\Biggl( \sum_{j=m}^k
a_j\mathbb{E} [V(X_j) ] \Biggr)^r\\
&\le& c_4^r \Biggl(\sum_{j=m}^k a_j j^{2\varepsilon} \Biggr)^r V^r(x_0)
\end{eqnarray*}
by Jensen's inequality and estimate (\ref{eq:drift-estimate}).
\end{pf}
%

Assume that $\{f_s\}_{s\in\mathbb{S}}$ is a regular enough family of
functions.
Consider the following decomposition, which is one of the key
observations in
\cite{andrieu-moulines},
%
%
\begin{equation} \label{eq:martingale-decomposition}
\sum_{j=1}^k [f_{S_j}(X_j) - \pi(f_{S_j}) ]
= M_k + R^{(1)}_k + R^{(2)}_k,
\end{equation}
where $(M_k)_{k\ge1}$ is a martingale with respect to $\F$,
and $(R^{(1)}_k)_{k\ge1}$ and
$(R^{(2)}_k)_{k\ge1}$ are ``residual'' sequences, given by
\begin{eqnarray*}
M_k &\defeq& \sum_{j=1}^k [
\hat{f}_{S_{j-1}}(X_j) - P_{S_{j-1}}\hat{f}_{S_{j-1}}(X_{j-1})
], \\
R^{(1)}_k &\defeq& \sum_{j=1}^k [
\hat{f}_{S_{j}}(X_j) - \hat{f}_{S_{j-1}}(X_{j})
], \\
R^{(2)}_k &\defeq&
P_{S_0}\hat{f}_{S_{0}}(X_0) - P_{S_k}\hat{f}_{S_{k}}(X_{k}).
\end{eqnarray*}
Recall that $\hat{f}_s$ solves the Poisson equation (\ref{eq:poisson}).
The following proposition controls the fluctuations of these terms
individually.
%
%
\begin{proposition}
\label{prop:main-andrieu} 
Assume \textup{\hyperlink{a:invariance}{(A1)}--\hyperlink{a:adapt-bound}{(A4)}} hold,
$(x_0,s_0)\in\mathbb{X}\times K_0$
and let $\{f_s\}_{s\in
\mathbb{S}}$ be $(K_n, V^\alpha)$-polynomially Lipschitz with constants
$(c,\varepsilon)$ for some $\alpha\in(0,1-\beta)$.
Then, for any $p\in(1,(\alpha+\beta)^{-1}]$,
for all $\delta>0$
and $\xi>\alpha$, there is
a $c_* = c_*(c,p,\alpha,\beta,\xi)\ge1$, such that for all $n\ge1$,
%
%
\begin{eqnarray}
\label{eq:main-thm-m}
\P\biggl[ \sup_{k\ge n} \frac{| M_k|}{k} \ge\delta\biggr]
&\le& c_* \delta^{-p}
n^{p\varepsilon_*-(p/2)\wedge(p-1)} V^{\alpha p}(x_0),
\\
\label{eq:main-thm-r1}
\P\biggl[ \sup_{k\ge n} \frac{|R_k^{(1)}|}{k^\xi}\ge\delta\biggr]
&\le& c_* \delta^{-p}
\Biggl(\sum_{j=1}^\infty
(j\vee n)^{\varepsilon_*-\xi}\eta_j \Biggr)^p
V^{(\alpha+\beta)p}(x_0),
\\
\label{eq:main-thm-r2}
\P\biggl[ \sup_{k\ge n} \frac{|R_k^{(2)}|}{k^\xi}\ge\delta\biggr]
&\le& c_* \delta^{-p}
n^{p\varepsilon_*-(\xi-\alpha)p} V^{\alpha p}(x_0),
\end{eqnarray}
whenever $\varepsilon>0$ is small enough to ensure that
$\varepsilon_* \defeq\kappa_* \varepsilon<
[\frac{1}{2} \wedge(1- \frac{1}{p} ) \wedge(\xi-\alpha) ]$,
where $\kappa_*\ge1$ is an independent constant.
\end{proposition}
%
%
\begin{pf} 
In this proof, $\tilde{c}$ is a constant that can take different
values at
each appearance.
By Proposition \ref{prop:poisson-regularity}, we have that
$\|\hat{f}_s\|_{V^\alpha} + \| P_s\hat{f}_s\|_{V^\alpha} \le c_3
\ell^{\kappa_3\varepsilon}$ for all $s\in K_\ell$. Since $\alpha p
\in
[0,1]$, we can bound
the martingale differences
$\ud M_\ell\defeq M_\ell- M_{\ell-1}$ for $\ell\ge1$
as follows:
%
%
\begin{eqnarray} \label{eq:dm-estimate}
\E
|\ud M_\ell|^p
&=& \E| \hat{f}_{S_{\ell-1}}(X_\ell)
- P_{S_{\ell-1}}\hat{f}_{S_{\ell-1}}(X_{\ell-1}) |^p\nonumber\\
&\le&\E\bigl|
\|\hat{f}_{S_{\ell-1}}\|_{V^\alpha} V^\alpha(X_\ell)
+ \| P_{S_{\ell-1}}\hat{f}_{S_{\ell-1}}\|_{V^\alpha} V^\alpha
(X_{\ell-1})
\bigr|^p\nonumber\\[-8pt]\\[-8pt]
&\le& 2^p(c_3 \ell^{\kappa_3\varepsilon})^p \bigl(
\mathbb{E} [V^{\alpha p}(X_\ell) ] + \mathbb{E} [V^{\alpha
p}(X_{\ell-1}) ]
\bigr)\nonumber\\
&\le& 2^{p+1}c_3^p c_4^{\alpha p}
\ell^{p\kappa_3\varepsilon}\ell^{2\alpha p\varepsilon} V^{\alpha p}(x_0)
\le\tilde{c} \ell^{(\kappa_3 + 2\alpha)p\varepsilon} V^{\alpha
p}(x_0)\nonumber
\end{eqnarray}
by (\ref{eq:stability1}) of Lemma \ref{lemma:stability}. For $p\ge2$,
we have, by Burkholder and Minkowski's inequalities,
\begin{eqnarray*}
\E|M_k|^p
&\le& c_p \mathbb{E} \Biggl[\sum_{\ell=1}^k |\ud M_\ell|^2 \Biggr]^{p/2}
\le
c_p \Biggl[\sum_{\ell=1}^k (\E
|\ud M_\ell|^p )^{2/p} \Biggr]^{p/2}\\
&\le&
\tilde{c} k^{(\kappa_3 + 2\alpha)p\varepsilon+ p/2} V^{\alpha p}(x_0),
\end{eqnarray*}
where the constant $c_p$ depends only on $p$.
For $1 < p \le2$, the estimate (\ref{eq:dm-estimate}) yields, by
Burkholder's inequality,
\begin{eqnarray*}
\E|M_k|^p
&\le& c_p
\mathbb{E} \Biggl[\sum_{\ell=1}^k (|\ud M_\ell|^p)^{2/p} \Biggr]^{p/2}
\le c_p \sum_{\ell=1}^k \E|\ud M_\ell|^p\\
&\le&\tilde{c} k^{(\kappa_3 +2\alpha)p\varepsilon+ 1} V^{\alpha p}(x_0).
\end{eqnarray*}
The two cases combined give that
%
%
\begin{equation} \label{eq:m-ex-p}
\E|M_k|^p
\le\tilde{c} k^{(\kappa_3+2\alpha)p\varepsilon+ (p/2)\vee1}
V^{\alpha
p}(x_0).
\end{equation}
Now, by Corollary
\ref{cor:birnbaum-marshall-martingale} of Birnbaum and Marshall's inequality
in Appendix~\ref{sec:martingale},
\begin{eqnarray*}
\P\biggl[
\max_{n\le k \le m} \frac{|M_k|}{k} \ge\delta
\biggr]
&\le&\delta^{-p} \Biggl[
m^{-p} \E|M_m|^p +
\sum_{k=n}^{m-1} \bigl(k^{-p} - (k+1)^{-p} \bigr) \E|M_k|^p
\Biggr] \\
&\le&\delta^{-p} \Biggl[
m^{-p} \E|M_m|^p +
p\sum_{k=n}^{m-1} k^{-p-1} \E|M_k|^p
\Biggr]
\end{eqnarray*}
for all $m\ge n$. By letting $\kappa_* \defeq\kappa_3 + 3$, we have
from (\ref{eq:m-ex-p})
\[
m^{-p} \E|M_m|^p
\le\tilde{c} m^{p(\kappa_*\varepsilon+ (1/2)\vee(1/p) - 1)}
\stackrel{m\to\infty}{\hbox to 1cm{\rightarrowfill}}0,
\]
since $\kappa_*\varepsilon+ (1/2)\vee(1/p) < 1$.
Now, (\ref{eq:main-thm-m}) follows by
\begin{eqnarray*}
\P\biggl[
\sup_{k\ge n} \frac{|M_k|}{k} \ge\delta
\biggr]
&\le&\tilde{c} \delta^{-p} \Biggl[
\sum_{k=n}^{\infty}
k^{(\kappa_3+2\alpha)p\varepsilon+ (p/2)\vee1 - p - 1}
\Biggr]
V^{\alpha p}(x_0)\\
&\le&\tilde{c} \delta^{-p}
n^{p\kappa_*\varepsilon- (p/2)\wedge(p - 1) }
V^{\alpha p}(x_0)
\end{eqnarray*}
since we have that
$p\kappa_*\varepsilon- (p/2)\wedge(p - 1)<0$.

By Proposition \ref{prop:poisson-regularity},
$\|\hat{f}_s-\hat{f}_{s'}\|_{V^\alpha} \le c_3
\ell^{\kappa_3\varepsilon}|s -s'|$
for $s,s'\in K_\ell$. By construction, $|S_\ell-S_{\ell-1}|\le\eta
_{\ell}
|H(S_{\ell-1},X_\ell)|$,
and assumption
\hyperlink{a:adapt-bound}{(A4)} ensures that
$|H(S_{\ell-1},X_\ell)|\le c \ell^\varepsilon V^\beta(X_\ell)$, so
\[
|\hat{f}_{S_\ell}(X_\ell)-\hat{f}_{S_{\ell-1}}(X_\ell) |
\le c_3 \ell^{\kappa_3\varepsilon}|S_\ell-S_{\ell-1}| V^\alpha
(X_\ell)
\le c_3 \ell^{\kappa_3\varepsilon}\eta_\ell c \ell^{\varepsilon}
V^{\alpha+\beta}(X_\ell).
\]
Let $k\ge n$. Since $\ell^{(\kappa_3+1)\varepsilon}k^{-\xi} \le
(\ell\vee n)^{(\kappa_3+1)\varepsilon-\xi}$ for $\ell\le k$, we obtain
\begin{eqnarray*}
k^{-\xi} \bigl|R_k^{(1)}\bigr|
&\le& k^{-\xi} \sum_{\ell=1}^k | \hat{f}_{S_\ell}(X_\ell)-\hat
{f}_{S_{\ell-1}}(X_\ell) |
\\
&\le&\tilde{c} \sum_{\ell=1}^k (\ell\vee n)^{(\kappa
_3+1)\varepsilon-\xi
} \eta_\ell
V^{\alpha+\beta}(X_\ell)
\end{eqnarray*}
and then by Minkowski's inequality and (\ref{eq:stability1}) of Lemma
\ref{lemma:stability},
%
%
\begin{eqnarray} \label{eq:exp-estim-r1}
&&
\E \Bigl[\max_{n\le k\le m} k^{-\xi p} \bigl|R^{(1)}_k\bigr|^p\Bigr]
\nonumber\\
&&\qquad\le\E\Biggl[
\sum_{\ell=1}^m \tilde{c} (\ell\vee n)^{(\kappa_3+1)\varepsilon
-\xi}
\eta_\ell
V^{(\alpha+\beta)p}(X_\ell)
\Biggr]^p \nonumber\\[-8pt]\\[-8pt]
&&\qquad\le\tilde{c} \Biggl[ \sum_{\ell=1}^m \bigl(\E
\bigl[(\ell\vee n)^{(\kappa_3+1)\varepsilon-\xi}
\eta_\ell
V^{\alpha+\beta}(X_\ell) \bigr]^p \bigr)^{1/p} \Biggr]^p \nonumber\\
&&\qquad\le\tilde{c}
\Biggl[\sum_{\ell=1}^\infty
(\ell\vee n)^{(\kappa_3+1+2\alpha+2\beta)\varepsilon-\xi} \eta
_\ell\Biggr]^p
V^{(\alpha+\beta)p}(x_0).\nonumber
\end{eqnarray}

Finally, consider $R_k^{(2)}$. From Proposition
\ref{prop:poisson-regularity}, we have that
$\|P_{S_k}\hat{f}_{S_k}(X_k)\|_{V^\alpha} \le
c_3 k^{\kappa_3\varepsilon}$,
and by (\ref{eq:stability2}) of Lemma \ref{lemma:stability},
\begin{eqnarray*}
\E\Bigl[
\max_{n\le k\le m} k^{-\xi p} |P_{S_k}\hat{f}_{S_k}(X_k)|^p
\Bigr]
&\le& c_3^p \E\Bigl[ \max_{n\le k\le m}
\bigl(k^{(\kappa_3\varepsilon-\xi)/\alpha} V(X_k) \bigr)^{\alpha p}
\Bigr] \\
&\le& c_3^{p}c_4^{\alpha p}
\Biggl(
\sum_{k=n}^m k^{(\kappa_3\varepsilon-\xi)/\alpha+2\varepsilon}
\Biggr)^{\alpha p}
V^{\alpha p}(x_0)
\\
&\le& \tilde{c} n^{ (\kappa_3+2\alpha)p\varepsilon+(\alpha-\xi) p}
V^{\alpha p}(x_0)
\end{eqnarray*}
since $(\kappa_3 + 2\alpha)\varepsilon- (\xi-\alpha) < 0$.
So, we have that
%
%
\begin{eqnarray} \label{eq:exp-estim-r2}\quad
\mathbb{E} \Bigl[\sup_{k\ge n} k^{-\xi p} \bigl| R^{(2)}_k \bigr|^p \Bigr]
&\le& 2^p
\mathbb{E} \Bigl[\sup_{k\ge n} k^{-\xi p} \bigl( |P_{S_0}\hat{f}_{S_0}(X_0)|^p
+ |P_{S_k}\hat{f}_{S_k}(X_k)|^p \bigr) \Bigr]
\nonumber\\
&\le&
2^p
\mathbb{E} \Bigl[|P_{S_0}\hat{f}_{S_0}(X_0)|^p + \sup_{k\ge n} k^{-\xi p}
|P_{S_k}\hat{f}_{S_k}(X_k)|^p \Bigr] \\
&\le& \tilde{c} n^{ (\kappa_3+2\alpha)p\varepsilon+(\alpha-\xi) p}
V^{\alpha p}(x_0).\nonumber
\end{eqnarray}
The estimates
(\ref{eq:main-thm-r1}) and (\ref{eq:main-thm-r2}) follow
by Markov's inequality from
(\ref{eq:exp-estim-r1}) and
(\ref{eq:exp-estim-r2}).
\end{pf}
%

The proof of Theorem \ref{thm:slln-restricted} follows as a
straightforward application of Proposition~\ref{prop:main-andrieu}.
%
\begin{pf*}{Proof of Theorem \ref{thm:slln-restricted}} 
Let $\delta>0$, and denote
\[
{B}_n^{(\delta)} \defeq
\Biggl\{ \omega\in\Omega:
\sup_{k\ge n}
\frac{1}{k} \Biggl|
\sum_{j=1}^k [f(X_j) -\pi(f) ] \Biggr| \ge\delta
\Biggr\}.
\]
Since $\|f\|_{V^\alpha}<\infty$ by assumption, we may consider the
family $\{f_s\}_{s\in\mathbb{S}}$ with $f_s \equiv f$ for all $s\in
\mathbb{S}$.
Then, we have by decomposition (\ref{eq:martingale-decomposition}) that
%
%
\begin{equation} \label{eq:slln-estimation}\qquad
\P\bigl(B_n^{(\delta)} \bigr)
\le
\P\biggl[\sup_{k\ge n} \frac{| M_k |}{k} \ge\frac{\delta}{3} \biggr]
+\P\biggl[\sup_{k\ge n} \frac{| R_k^{(1)} |}{k} \ge\frac{\delta}{3} \biggr]
+\P\biggl[\sup_{k\ge n} \frac{| R_k^{(2)} |}{k} \ge\frac{\delta}{3} \biggr].
\end{equation}
We select $p\in(1,(\alpha+\beta)^{-1})$ so that
$\kappa_*\varepsilon<(1-1/p)$, and let $\xi=1$. Then, Proposition
\ref{prop:main-andrieu} readily implies that the first and the third
terms in
(\ref{eq:slln-estimation}) converge to zero as $n\to\infty$. For the
second term, consider
\[
\sum_{j=1}^\infty
(j\vee n)^{\kappa_*\varepsilon-1}\eta_j
= n^{\kappa_*\varepsilon-1}\sum_{j=1}^n \eta_j + \sum
_{j=n+1}^\infty
j^{\kappa_*\varepsilon-1}\eta_j,
\]
where the second term converges to zero by assumption, and the
first term by Kronecker's lemma.
There is an increasing sequence $(n_k)_{k\ge1}$ such that
$\P({B}_{n_k}^{(1/k)})\le k^{-2}$. Denoting ${B}\defeq
\bigcap_{m=1}^\infty\bigcup_{k=m}^\infty{B}_{n_k}^{(1/k)}$, the
Borel--Cantelli lemma implies that $P(B^\complement)=1$,
and for all $\omega\in{B}^\complement$, (\ref{eq:slln}) holds.
\end{pf*}

\section{Bound for the growth rate}
\label{sec:grow-bound} 

In this section, we assume that $\mathbb{X}$ is a normed space, and
establish a
bound for the growth rate of the chain $(\|X_n\|)_{n\ge1}$, based
on a
general drift condition. The bound assumes little structure; one must
have a
drift function $V$ that grows rapidly enough, and that the expected growth
of $V(X_n)$ is moderate.
%
\begin{proposition}
\label{prop:grow-bound} 
Suppose that there is $V:\mathbb{X}\to[1,\infty)$ such
that the bound
%
%
\begin{equation} \label{eq:drift-bound}
P_s V(x) \le V(x) + b
\end{equation}
holds for all $(x,s)\in\mathbb{X}\times\mathbb{S}$, where $b<\infty
$ is a
constant independent of $s$.
Suppose also that $V$ grows rapidly enough so that
%
%
\begin{equation} \label{eq:norm-v-bound-1}
\|x\| \ge u \quad\Longrightarrow\quad
V(x)\ge r(u)
\end{equation}
for all $u\ge0$, where $r\dvtx[0,\infty)\to[0,\infty)$ is a function
growing faster than any polynomial, that is, for any $p>0$ there is
a $c=c(p)<\infty$ such that
%
%
\begin{equation} \label{eq:norm-v-bound-2}
\sup_{u\ge1}\frac{u^p}{r(u)} \le c.
\end{equation}
Then, for any $\varepsilon>0$,
there is an a.s. finite $A=A(\omega,\varepsilon)$ such that
\[
\|X_n\| \le An^{\varepsilon}.
\]
\end{proposition}
%
%
\begin{pf} 
To start with, (\ref{eq:drift-bound}) implies for $n\ge1$
\begin{eqnarray*}
\mathbb{E} [V(X_n) ] &=& \mathbb{E} [\mathbb{E}[V(X_n) | \F_{n-1}] ]
= \mathbb{E} [P_{S_{n-1}}V(X_{n-1}) ]
\le\mathbb{E} [V(X_{n-1}) ] + b \\
&\le&\cdots\le V(x_0) + b n \le\tilde{b}V(x_0) n,
\end{eqnarray*}
where $\tilde{b} \defeq b+1$.
Now, with fixed $a\ge1$, we can bound the probability of
$\|X_n\|$ ever exceeding $an^\varepsilon$ as follows
\begin{eqnarray*}
\P\biggl(
\max_{1\le n\le m} \frac{\|X_n\|}{n^\varepsilon} \ge a
\biggr)
&\le&
\sum_{n=1}^m
\P(
\|X_n\| \ge an^{\varepsilon}
)
\le\sum_{n=1}^\infty
\P\bigl(
V(X_n) \ge r(an^{\varepsilon})
\bigr) \\
&\le&
\sum_{n=1}^\infty
\frac{\mathbb{E} [V(X_n) ]}{r(an^\varepsilon)}
\le
\tilde{b} V(x_0)
\sum_{n=1}^\infty
\frac{n}{r(an^\varepsilon)}\\
&\le&
\frac{\tilde{b}V(x_0)c}{a^{3/\varepsilon}}
\sum_{n=1}^\infty
n^{-2}
\stackrel{a\to\infty}{\hbox to 1cm{\rightarrowfill}}0,
\end{eqnarray*}
where we use Markov's inequality, and
$c=c(3/\varepsilon)<\infty$ is from the application of
(\ref{eq:norm-v-bound-2}).
\end{pf}
%

We record the following easy lemma, dealing with a particular choice of $V(x)$,
for later use in Section \ref{sec:am-ergodicity}.
%
\begin{lemma}
\label{lemma:v-polybound} 
Assume that the target density $\pi$ is differentiable, bounded, bounded
away from zero on compact sets, and satisfies the following radial
decay condition:
\[
\lim_{r\to\infty}\sup_{\|x\|\ge r}
\frac{x}{\|x\|} \cdot\nabla\log\pi(x) < 0.
\]
Then, for
$V(x)=c_V\pi^{-1/2}(x)$, the bound (\ref{eq:norm-v-bound-1}) applies
with a function $r(u) \defeq c e^{\gamma u}$ for some $\gamma,c>0$,
satisfying (\ref{eq:norm-v-bound-2}).
\end{lemma}
%
%
\begin{pf} 
Let $R\geq1$ be such that $\sup_{\|x\|\ge R}
\frac{x}{\|x\|} \cdot\nabla\log\pi(x)\le-\gamma$ for some
$\gamma>0$.
Assume $y\in\R^d$ and $\|y\|\ge2R$, and write
$y = (1+a)x$, where $\|x\| = R$ and $a=\frac{\|y\|}{R}-1\ge1$.
Denote $h(x)\defeq\log\pi(x)$, and write
\[
\log
\frac{\pi(y)}{\pi(x)} = \int_1^{1+a}
{x} \cdot\nabla h(t x)\,dt
\le-\gamma a.
\]
We have that
\[
V(y)
= c_V \pi(x)^{-1/2} \biggl(\frac{\pi(y)}{\pi(x)} \biggr)^{-1/2}
\ge c_V e^{{\gamma a}/{2}} \inf_{\|x\|=R} \pi(x)^{-1/2}
\ge c e^{{\gamma}/({4R})\|y\|}
\]
and, since $\pi$ is bounded away from zero on $\{x\dvtx\|x\|< 2R\}$,
we can select $c>0$ such that the bound applies to all $y\in\R^d$.
\end{pf}
%

\section{Ergodicity result for adaptive metropolis}
\label{sec:am-ergodicity} 

We start this section by outlining the original Adaptive Metropolis
(AM) algorithm
\cite{saksman-am}. The AM chain starts from a point $X_0 \equiv x_0\in
\R^d$,
and we have an initial covariance
$\Sigma_0\in\mathcal{C}^d$ where $\mathcal{C}^d\subset\R^{d\times d}$
stands for the symmetric and positive definite matrices.
We generate, recursively, for $n\ge0$,
%
%
\begin{eqnarray}
X_{n+1} &\sim& P_{\theta\Sigma_n}(X_n, \uarg),
\\
\label{eq:am-orig-s}
\Sigma_{n+1} &=& \cases{v_0,&\quad $0\le n\le N_b-1$,\cr
\operatorname{Cov}(X_0,\ldots,X_n)+\kappa I,&\quad $n\ge
N_b$,}
\end{eqnarray}
where $\theta>0$ is a parameter,
$N_b\ge2$ is the length of the burn-in,
$\kappa>0$ is a small constant,
$I$ is an identity matrix
and $P_v(x,\uarg)$ is a Metropolis transition probability defined as
%
%
\begin{eqnarray} \label{eq:metropolis-kernel}
P_v(x,{A}) &\defeq&
\mathbh{1}_{{A}}(x) \biggl[ 1-
\int\biggl(1\wedge\frac{\pi(y)}{\pi(x)} \biggr) q_v(y-x) \,\ud y \biggr]
\nonumber\\[-8pt]\\[-8pt]
&&{}+ \int_{{A}}
\biggl(1\wedge\frac{\pi(y)}{\pi(x)} \biggr) q_v(y-x) \,\ud y,
\nonumber
\end{eqnarray}
where the proposal
density $q_v$ is the Gaussian density with zero mean and covariance
$v\in\mathcal{C}^d$.

In this paper, just for notational simplicity (see Remark \ref
{rem:convergence}),
we consider a slight modification of the AM chain.
First, we do not consider a burn-in period, that is, let $N_b=0$, and let
$\Sigma_0\ge\kappa I$.
Instead of (\ref{eq:am-orig-s}), we construct $\Sigma_{n}$
recursively for $n\ge1$ as
%
%
\begin{equation} \label{eq:mod-recursion}
\Sigma_{n} = \frac{n}{n+1} \Sigma_{n-1}
+ \frac{1}{n+1}
[(X_{n}-\overline{X}_{n-1})(X_{n}-\overline{X}_{n-1})^T
+ \kappa I ],
\end{equation}
where $\overline{X}_n$ denotes the average of $X_0,\ldots,X_n$.
%
\begin{remark}
\label{rem:convergence} 
The original AM process uses the unbiased estimate of the
covariance matrix. In this case, the recursion formula for $\Sigma_n$, when
$n\ge N_b+2$, has the form
%
%
\begin{equation} \label{eq:am-orig-crec}
\Sigma_n = \frac{n-1}{n} \Sigma_{n-1}
+ \frac{1}{n+1} [
(X_n - \overline{X}_{n-1})(X_n - \overline{X}_{n-1})^T + \kappa I ].
\end{equation}
This recursion can also be formulated in our framework described in
Section \ref{sec:notations} by simply introducing a sequence of adaptation
functions $H_n(s,x)$. Our proof applies with obvious changes. However, in
the present paper, we prefer (\ref{eq:mod-recursion}) for simpler notation.
Also, from a practical point of view, observe that (\ref
{eq:mod-recursion}) differs from
(\ref{eq:am-orig-crec}) by a factor smaller than $n^{-2} \Sigma_{n-1}$
whence it is mostly a matter of taste whether to use
(\ref{eq:mod-recursion}) or (\ref{eq:am-orig-crec}).
\end{remark}
%

In the notation of the general adaptive MCMC framework in
Section \ref{sec:notations}, we have the state space $\mathbb
{X}\defeq\R
^d$. The
adaptation parameter $S_n=(S_n^{(m)},S_n^{(v)})$ consists of the mean
$S_n^{(m)}$ and the covariance $S_n^{(v)}$, having values in
$(S_n^{(m)},S_n^{(v)}) \in\mathbb{S} \defeq\R^d\times\mathcal{C}^d$.
The space $\overline{\mathbb{S}}\defeq\R^d\times\R^{d\times d}
\supset
\mathbb{S}$
is equipped with the norm $|s| \defeq\|s^{(m)}\| \vee
\|s^{(v)}\|$ where we use the Euclidean norm, and the matrix norm
$\|A\|^2\defeq\trace(A^T A)$, respectively.
The Metropolis kernel $P_s$ is defined as in (\ref{eq:metropolis-kernel}),
with the definition $q_s \defeq q_{s^{(v)}}$ for
$s\in\mathbb{S}$.
The adaptation function $H$ is
defined for $s=(s^{(m)},s^{(v)})$ as
\[
H(s, x) \defeq\left[\matrix{
x - s^{(m)} \cr
\bigl(x-s^{(m)}\bigr)\bigl(x-s^{(m)}\bigr)^T - s^{(v)} + \kappa I}\right],
\]
and the adaptation weights are $\eta_n \defeq(n+1)^{-1}$.

We now formulate our ergodicity result for the AM chain.
%
\begin{theorem}
\label{thm:am-ergodicity} 
Assume $\pi$ is positive, bounded, bounded from
below on compact sets, differentiable and
%
%
\begin{equation} \label{eq:super-pexp}
\lim_{r\to\infty} \sup_{\|x\|\ge r}
\frac{x}{\|x\|^\rho} \cdot\nabla\log\pi(x) = -\infty
\end{equation}
for some $\rho>1$. Moreover, assume that $\pi$ has regular contours
%
%
\begin{equation} \label{eq:reg-contour}
\lim_{r\to\infty}
\sup_{\|x\|\ge r}
\frac{x}{\|x\|} \cdot\frac{\nabla\pi(x)}{\|\nabla\pi (x)\|} < 0.
\end{equation}

Define $V(x)\defeq c_V \pi^{-1/2}(x)$ with $c_V = (\sup_x \pi(x))^{1/2}$.
Then, for any $f$ with $\|f\|_{V^\alpha} < \infty$ where
$0\le\alpha< 1$,
%
%
\begin{equation} \label{eq:am-slln}
\frac{1}{n}\sum_{k=1}^n f(X_k) \stackrel{n\to\infty}{\hbox to
1cm{\rightarrowfill}} \pi(f)
\end{equation}
almost surely.
\end{theorem}
%
%
\begin{remark}
\label{rem:slln-moments} 
If the conditions of Theorem \ref{thm:am-ergodicity} are
satisfied, the function $V(x)$ grows faster than an exponential, and
hence (\ref{eq:am-slln}) holds for exponential moments.
In particular, (\ref{eq:am-slln}) holds for power moments, that is,
for $f(x) = \|x\|^p$
for any $p\ge0$, and therefore also $S_n \to(m_\pi, v_\pi+\kappa I)$
where $m_\pi$ and $v_\pi$ are the mean and covariance of $\pi$.
\end{remark}
%

The proof of Theorem \ref{thm:am-ergodicity} is postponed to the end
of this
section. We start by a simple lemma bounding the growth rate of the AM chain.
%
\begin{lemma}
\label{lemma:am-bound} 
If the conditions of Proposition \ref{prop:grow-bound}
are satisfied for an AM chain, then for any $\varepsilon>0$,
there is an a.s. finite $A=A(\omega,\varepsilon)$ such that
\[
\bigl\| S_n^{(m)}\bigr\| \le An^\varepsilon,\qquad
\bigl\| S_n^{(v)}\bigr\| \le An^{\varepsilon}.
\]
\end{lemma}
%
%
\begin{pf} 
Since the AM recursion is a convex combination, this
is a straightforward
corollary of Proposition \ref{prop:grow-bound}.
\end{pf}
%

Next, we show that each of the Metropolis kernels used by the AM algorithm
satisfy a geometric drift condition, and bound the constants of geometric
drift. The result in Proposition \ref{prop:geom-bound} is similar to the
results obtained in \cite{jarner-hansen,roberts-tweedie}, with the exception
that we have a common minorization set $C$ for all proposal scalings.
We start by two lemmas. We define $\ball(x,r) \defeq\{y\in\R^d\dvtx
\|x-y\|\le r\}$.
%
\begin{lemma}
\label{lemma:boundary-estim} 
Assume $E\subset\R^d$ is measurable and $A\subset\R^d$ compact,
given as
\[
A \defeq\{ ru \dvtx u\in S^d, 0\le r\le g(u) \},
\]
where $S^d\defeq\{u\in\R^d\dvtx\|u\|=1\}$ is the unit sphere, and
$g\dvtx S^d\to[b,\infty)$ is a measurable function parameterising
the boundary $\partial A$,
with some $b>0$.

For any $\varepsilon>0$, define $B_\varepsilon\defeq\{ru\dvtx u\in S^d,
g(u)<r\le
g(u)+\varepsilon\}$. Then,
for all $\tilde{\varepsilon}>0$, there is
a $\tilde{b}=\tilde{b}(\tilde{\varepsilon})\in(0,\infty)$
such that for all $0<\varepsilon<\tilde{\varepsilon}$ and
for all $\lambda\ge3 \varepsilon$,
it holds that
\[
|E \cap B_\varepsilon| \le
\bigl| \bigl(E \oplus\ball(0,\lambda) \bigr)\cap A \bigr|,
\]
whenever $b\ge\tilde{b}$. Above, $A\oplus B \defeq\{x+y\dvtx x\in A,
y\in
B\}$ stands for the Minkowski sum.
\end{lemma}
%
%
\begin{pf} 
See Figure \ref{fig:contour-est} for an illustration of the situation.
Denote by $S^* \defeq\{u\in S^d\dvtx\exists r>0, ur\in E\cap
B_\varepsilon\}$
the projection of the set $E\cap B_\varepsilon$ onto $S^d$.
%
%
\begin{figure}

\includegraphics{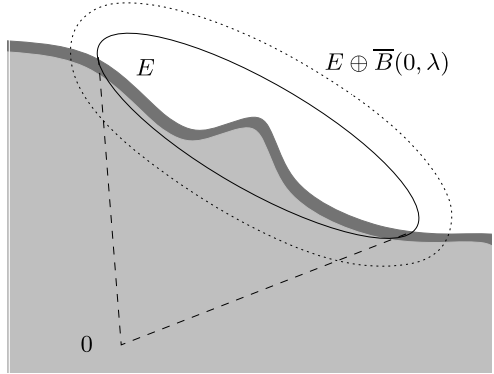}

\caption{Illustration of the boundary estimate. The set $A$ is
in light grey, and the set $B_\varepsilon$ in dark gray.} 
\label{fig:contour-est}
\end{figure}
Then we have $E\cap B_{\varepsilon} \subset\{ru\dvtx u\in S^*,
g(u)<r\le
g(u)+\varepsilon\}$
and $A\supset\{ru\dvtx u\in S^*, 0\le r\le
g(u)\}$. Now, for $\varepsilon\le\lambda\le g(u)$, we have
\[
\bigl((E\cap B_\varepsilon)\oplus\ball(0,\lambda) \bigr)\cap A
\supset\{ru\dvtx u\in
S^*, g(u)-\lambda+\varepsilon\le r\le g(u)\} \eqdef G,
\]
for let $ru \in G$, then there is $g(u)< \tilde{r} \le g(u) +
\varepsilon$ such
that $\tilde{r}u \in E\cap B_\varepsilon$, and we can write $ru =
\tilde
{r}u +
(r-\tilde{r})u$, where $(r-\tilde{r})u\in\ball(0,\lambda)$.
Clearly, $E \oplus\ball(0,\lambda)
\supset(E \cap B_\varepsilon) \oplus\ball(0,\lambda)$, and we
can estimate
\begin{eqnarray*}
&&
\bigl| \bigl(E \oplus\ball(0,\lambda) \bigr)\cap A \bigr|
- |E \cap B_\varepsilon| \\
&&\qquad\ge\int_{S^*} \int_{g(u)-2\varepsilon}^{g(u)} r^{d-1} \,\ud r
- \int_{g(u)}^{g(u)+\varepsilon}r^{d-1} \,\ud r\, \mathcal{H}^{d-1}(\ud u)
\\
&&\qquad=\frac{1}{d}\int_{S^*}
2 (g(u))^d - \bigl(g(u)-2\varepsilon\bigr)^d - \bigl(g(u)+\varepsilon\bigr)^d
\mathcal{H}^{d-1}(\ud u),
\end{eqnarray*}
where $\mathcal{H}^{d-1}$ stands for the $(d-1)$-dimensional Hausdorff
measure.
This integral is nonnegative for all $0\le\varepsilon\le c_d b$, for some
constant $c_d$ depending only on the dimension $d$,
namely let $h(\varepsilon)\defeq(y-2\varepsilon)^d + (y+\varepsilon
)^d$. The mean
value theorem implies that for some $0\le\varepsilon'\le\varepsilon$,
one has
\[
h(0)-h(\varepsilon)
= \varepsilon d (y-2\varepsilon')^{d-1} \biggl[ 2 -
\biggl(\frac{y+\varepsilon'}{y-2\varepsilon'} \biggr)^{d-1} \biggr]
\ge0,
\]
whenever $\varepsilon\le c_d y$.
\end{pf}
%
%
%
\begin{lemma}
\label{lemma:exp-estim} 
Let $f(x) \defeq x e^{-{x^2}/{2}}$.
For any $0<\varepsilon<1/8$, the following estimates hold:
\[
2f(x+\varepsilon) - f(x) \ge \frac{x}{8}\qquad \mbox{for all } 0<x\le
\frac{1}{2}
\]
and
\[
\int_0^\infty\bigl([2f(x+\varepsilon) - f(x)]\wedge0 \bigr) \,\ud x
\ge - e^{-c \varepsilon^{-2}}
\]
for some constant $c>0$.
\end{lemma}
%
%
\begin{pf} 
We can write
\[
2 f(x+\varepsilon) - f(x) = e^{-{x^2}/{2}} [ 2(x+\varepsilon
)e^{-x\varepsilon
- {\varepsilon^2}/{2}} - x ],
\]
which is positive whenever $e^{-x\varepsilon
- {\varepsilon^2}/{2}}\ge2/3$, holding at least
for all $0\le x \le x^*$, with
\[
x^* = \frac{\log(3/2)}{\varepsilon} - \frac{\varepsilon}{2}
\ge\frac{1}{4\varepsilon}.
\]
Now, $x^*\ge1/2$ and we can estimate
\[
2 f(x+\varepsilon) - f(x)
\ge\frac{1}{4} x e^{-{x^2}/{2}} \ge\frac{x}{8}
\]
for all $0<x\le1/2$. Also,
\[
\int_0^\infty\bigl([2f(x+\varepsilon) - f(x)]\wedge0 \bigr) \,\ud x
\ge-\int_{x^*}^\infty x e^{-{x^2}/{2}} \,\ud x
= -e^{-c\varepsilon^{-2}}
\]
with $c = 1/32$.
\end{pf}
%
%
%
\begin{proposition}
\label{prop:geom-bound} 
Assume that $\pi$ satisfies the conditions in Theorem
\ref{thm:am-ergodicity} and $\kappa>0$.
Then, there exists a compact set
$C\subset\R^d$, a probability measure $\nu$ on $C$,
and a constant
$b\in[0,\infty)$ such that
for the Metropolis transition probability $P_v$
in (\ref{eq:metropolis-kernel}) and
for all $v\in\mathcal{C}^d$ with all
eigenvalues greater than $\kappa>0$, it holds that
%
%
\begin{eqnarray}
\label{eq:drift-ineq-m}
P_v V(x) &\le& \lambda_v V(x) + b \mathbh{1}_{C}(x)\qquad
\forall x\in\mathbb{X},\\
\label{eq:mino-ineq-m}
P_v(x,B) &\ge& \delta_v \nu({B}) \qquad \forall x\in C, \forall
{B}\subset\mathbb{X},
\end{eqnarray}
where
$V(x)\defeq c_V \pi^{-1/2}(x)\ge1$ with
$c_V \defeq(\sup_x \pi(x))^{1/2}$ and the constants
$\lambda_v,\delta_v\in(0,1)$ satisfy the bound
\[
(1-\lambda_v)^{-1} \vee\delta_v^{-1} \le c \det(v)^{1/2}
\]
for some constant $c\ge1$.
\end{proposition}
%
%
\begin{pf} 
Define the sets $A_x\defeq\{y\dvtx\pi(y)\ge\pi(x)\}$ and its complement
$R_x\defeq\{y\dvtx\pi(y)<
\pi(x)\}$, which are the regions of almost sure acceptance and possible
rejection at $x$, respectively.
Let $R>1$ be sufficiently large to ensure that for all $\|x\|\ge
R$, it
holds that
\[
\sup_{\|x\|\ge R}
\frac{x}{\|x\|} \cdot\frac{\nabla\pi(x)}{\|\nabla\pi (x)\|} <
-\gamma
\quad\mbox{and}\quad
\sup_{\|x\|\ge R} \frac{x}{\|x\|} \cdot\nabla
\log\pi(x) < -\|x\|^{\rho-1}
\]
for some $\gamma>0$.
Suppose that the dimension $d\ge2$.
Lemma \ref{lemma:contour} in Appendix \ref{sec:contour} implies that
for $R$ sufficiently large, we have
$\ball(0,M^{-1}\|x\|)\subset A_x\subset\ball(0,M\|x\|)$ for all
$\|x\|\ge R$
with some constant $M\ge1$.
Moreover, we can parameterize $A_x = \{ru\dvtx u\in S^d, 0\le r\le
g(u)\}$
where $S^d \defeq\{u\in\R^d\dvtx\|u\|=1\}$ is the unit sphere, and
$g\dvtx S^d\to[M^{-1}\|x\|,M\|x\|]$.

Consider (\ref{eq:drift-ineq-m}). We may compute
%
%
\begin{eqnarray}\label{eq:drift-fraction}
\tau_v :\!&=&
1-
\frac{P_v V(x)}{V(x)}\nonumber\\
&=& \int_{A_x} \Biggl( 1- \sqrt{\frac{\pi(x)}{\pi(y)}}
\Biggr) q_v(y-x) \,\ud y \\
&&{}- \int_{R_x} \sqrt{\frac{\pi(y)}{\pi(x)}}
\Biggl( 1 -
\sqrt{\frac{\pi(y)}{\pi(x)}} \Biggr) q_v(y-x) \,\ud y.\nonumber
\end{eqnarray}
In what follows, unless explicitly stated, we assume $\|x\|\ge M(R
+ 1)$.
Denote $\varepsilon_x \defeq\|x\|^{-\alpha} < 1$, where $\alpha=
(\rho-1)/2>0$.
Define $\tilde{A}_x \defeq\{ru\dvtx u\in S^d, 0\le r\le
g(u)-\varepsilon_x\}\subset A_x$ and $\tilde{R}_x \defeq\{ru\dvtx
u\in
S^d, r\ge g(u)+\varepsilon_x \}\subset R_x$.
From (\ref{eq:drift-fraction}), we can estimate
%
%
\begin{eqnarray}\label{eq:tau-estim}
\tau_v &\ge&
\int\Biggl[
\Biggl( 1- \sqrt{\frac{\pi(x)}{\pi(y)}} \Biggr) \mathbh{1}_{\tilde{A}_x}(y)
- \frac{1}{4} \mathbh{1}_{R_x\setminus\tilde{R}_x}(y)
\Biggr]q_v(y-x) \,\ud y \nonumber\\[-8pt]\\[-8pt]
&&{}- \sup_{z\in\R^d} q_v(z-x) \int_{\tilde{R}_x} \sqrt{\frac{\pi
(y)}{\pi(x)}}
\,\ud y.\nonumber
\end{eqnarray}
We estimate the two terms in the right-hand side separately,
starting from the first.

Let $h(x)\defeq\log\pi(x)$.
Suppose $z\in\tilde{A}_x$, and write $z = (1-a/\|y\|)y$
for some $y\in\partial A_x$ and $\varepsilon_x \le a \le\|y\|$.
Assume for a moment
$\|z\|\ge R$. Then, $h$ is decreasing on the line segment from $z$ to
$y$, and we can estimate
\begin{eqnarray*}
\frac{\pi(x)}{\pi(z)} &=&
\frac{\pi(y)}{\pi(z)} = e^{h(y)-h(z)}
= e^{\int_{\|y\|-a}^{\|y\|}
{y}/{\|y\|}\cdot\nabla h(t{y}/{\|y\|})\,\ud t}
\le e^{\int_{\|y\|-\varepsilon_x}^{\|y\|}
{y}/{\|y\|}\cdot\nabla h(t{y}/{\|y\|})\,\ud t} \\
&\le& e^{-\varepsilon_x(\|y\|-\varepsilon_x)^{\rho-1}}
\le e^{-\varepsilon_x\|x\|^{\rho-1}/(2M)^{\rho-1}}
= e^{-\|x\|^\alpha/(2M)^{\rho-1}}.
\end{eqnarray*}
Hence, in this case, $\pi(x)/\pi(z)\le1/4$ assuming $\|x\|\ge R_2$
for sufficiently large $R_2\ge R$.
If $\|z\|<R$, then there is $z'$ such that $\|z'\|=R$ and
the estimate above holds for $z'$. Consequently,
%
%
\begin{equation} \label{eq:ax-estimate}
\frac{\pi(x)}{\pi(z)} = \frac{\pi(y)}{\pi(z')}\frac{\pi
(z')}{\pi(z)}
\le e^{-\|x\|^\alpha/(2M)^{\rho-1}} \frac{\sup_{\|w\|\le
R}\pi(w)}{\inf_{\|w\|\le
R}\pi(w)}
\le\frac{1}{4},
\end{equation}
whenever $\|x\| \ge R_2$ by increasing $R_2$ if needed.
In conclusion, we have shown that for
$\|x\|\ge R_2$, it holds that
$(1-\sqrt{\pi(x)/\pi(y)}) \ge1/2$ for all $y\in\tilde{A}_x$.

By Fubini's theorem, we can write for positive $f$ that
\begin{eqnarray*}
\int f(z+x) q_v(z) \,\ud x
&=& \frac{c_d}{\sqrt{\det(v)}} \int_0^1
\int_{\{e^{-{1/2}z^Tv^{-1}z}\ge t\}}
f(z+x) \,\ud z \,\ud t \\
&=& \frac{c_d}{\sqrt{\det(v)}} \int_0^\infty
\int_{E_u}
f(y) \,\ud y\, u e^{-{u^2}/{2}}\,\ud u,
\end{eqnarray*}
where $c_d=(2\pi)^{-d/2}$ and
$E_u \defeq\{z+x\dvtx z^Tv^{-1}z\le u^2\}$.
Consequently, for $\|x\|\ge R_2$, we
can estimate the first term of (\ref{eq:tau-estim}) from below by
\begin{eqnarray*}
&&\int_0^\infty
\biggl(\frac{|E_u\cap\tilde{A}_x|}{2} - \frac{|E_u\cap(R_x\setminus
\tilde{R}_x)|}{4} \biggr) u e^{-{u^2}/{2}}\,\ud u
\\
&&\qquad\ge \frac{1}{4}\int_0^\infty
2|E_{u+a}\cap\tilde{A}_x| (u+a) e^{-{(u+a)^2}/{2}}
- |E_u\cap(R_x\setminus
\tilde{R}_x)| u e^{-{u^2}/{2}}\,\ud u \\
&&\qquad\ge \frac{1}{4}\int_0^\infty
2 \bigl| \bigl(E_{u}\oplus\ball(0,\kappa^{1/2}a) \bigr)\cap\tilde{A}_x \bigr|
(u+a) e^{-{(u+a)^2}/{2}}\\
&&\qquad\quad\hspace*{26.4pt}{}- |E_u\cap B_\varepsilon| u e^{-{u^2}/{2}}\,\ud u
\end{eqnarray*}
for any $a\ge0$, since simple computation shows that
$E_u \oplus\ball(0,\kappa^{1/2}a)
= \{x+y\dvtx x\in E_u, y\in\ball(0,\kappa^{1/2}a)\} \subset E_{u + a}$,
and as we may write
$\tilde{A}_x = \{ru\dvtx u\in S^d, 0\le r\le\tilde{g}(u)\}$
where $\tilde{g}(u) = g(u)-\varepsilon_x$, we obtain that
$R_x\setminus
\tilde{R}_x \subset\{ru\dvtx u\in S^d, \tilde{g}(u)\le r\le\tilde
{g}(u) +
2\varepsilon_x\} \eqdef B_\varepsilon$.
We set $a=6\kappa^{-1/2}\varepsilon_x$ and apply
Lemma \ref{lemma:boundary-estim} with the choice $\varepsilon
=2\varepsilon_x$ and
$\lambda=6\varepsilon_x$,
\begin{eqnarray*}
&&\int_0^\infty
\biggl(\frac{|E_u\cap\tilde{A}_x|}{2} - \frac{|E_u\cap(R_x\setminus
\tilde{R}_x)|}{4} \biggr) u e^{-{u^2}/{2}}\,\ud u
\\
&&\qquad \ge\frac{1}{4}\int_0^\infty
| [E_u\oplus\ball(0,6\varepsilon_x) ]\cap\tilde{A}_x |
\bigl[2 (u+a) e^{-{(u+a)^2}/{2}}
- u e^{-{u^2}/{2}} \bigr]\,\ud u \\
&&\qquad\ge\frac{1}{4}\int_{1/4}^{1/2}
|E_{u}\cap\tilde{A}_x| \frac{u}{8}\,\ud u - |\tilde{A}_x|
e^{-c_1\varepsilon_x^{-2}}
\\
&&\qquad\ge c_2 |E_{1/4}\cap\tilde{A}_x| - M^d\|x\|^d e^{-c_1\|x\|^{\alpha}}
\end{eqnarray*}
by Lemma \ref{lemma:exp-estim},
for sufficiently large $\|x\|$, and since $E_u$ are increasing with
respect to $u$.
We have that $E_{1/4}\supset\ball(x,\kappa^{1/2}/4)$.
If $\|x\|\to\infty$, then $\varepsilon_x\to0$ and also
$|\ball(x,\kappa^{1/2}/4)\cap\tilde{A}_x| -|\ball(x,\kappa
^{1/2}/4)\cap
A_x| \to0$. Moreover, it holds that
$|\ball(x,\kappa^{1/2}/4)\cap A_x|\ge c_3>0$ (see
the proof of Theorem 4.3 in \cite{jarner-hansen}).
So, for large enough $\|x\|$, there is a $c_4>0$ so that
$|E_{1/4}\cap\tilde{A}_x|\ge c_4$.
To sum up, by choosing $R_3$ to be sufficiently large, we obtain that the
first part of (\ref{eq:tau-estim}) is at least $c_5 (\det(v))^{-1/2}$ for
all $\|x\|\ge R_3$, with a $c_5>0$.

Next, we turn to the second term of (\ref{eq:tau-estim}).
We obtain by polar integration that
\begin{eqnarray*}
\int_{\tilde{R}_x} \sqrt{\frac{\pi(y)}{\pi(x)}}
\,\ud y
&=&\int_{S^d}
\int_{g(u)+\varepsilon_x}^\infty r^{d-1} e^{{1}/{2}h(r u)-
{1}/{2}h(g(u)u)}\, \ud r\,
\mathcal{H}^{d-1}(\ud u) \\
&\le& c_d' \sup_{M^{-1}\|x\|\le w\le M\|x\|}
\int_{w+\varepsilon_x}^\infty r^{d-1} e^{-{1}/{2}\int_{w}^r
t^{\rho-1} \,dt}\,\ud r,
\end{eqnarray*}
where $\mathcal{H}^{d-1}$ is the $(d-1)$-dimensional Hausdorff
measure, and
$c_d'=\mathcal{H}^{d-1}(S^d)$.
Denote $T(w,r) \defeq r^{d-1}
e^{-{1}/{4}\int_{w}^r t^{\rho-1}\,\ud t}$
and let us estimate the latter integral
from above by
\begin{eqnarray*}
\int_{w+\varepsilon_x}^\infty
e^{-{1}/{4}\int_{w}^r t^{\rho-1}\,\ud t} \,\ud r
\sup_{r\ge w+\varepsilon_x} T(w,r)
&\le&
\int_{w}^\infty e^{-{w^{\rho-1}}/{4}(r-w)} \,\ud r
\sup_{r\ge w+\varepsilon_x} T(w,r) \\
&\le&
4M^{\rho-1}\|x\|^{1-\rho}
\sup_{r\ge w+\varepsilon_x} T(w,r)
\end{eqnarray*}
for any $w\ge M^{-1}\|x\|$.
Suppose first $w+\varepsilon_x \le r \le2w$, then
\[
T(w,r) \le(2w)^{d-1} e^{-{1}/{4}\varepsilon_x w^{\rho-1}}
\le(2M)^{d-1} \|x\|^{d-1} e^{-{1}/{4}M^{1-\rho} \|x\|^\alpha}
\le c_6
\]
for any $M^{-1}\|x\| \le w \le M\|x\|$.
For any $r>2w$ and $w\ge1$, we have
\[
T(w,r) \le r^{d-1} e^{-{1}/{4} {r}/{2} w^{\rho-1}}
\le r^{d-1} e^{-{r}/{8}} \le c_7.
\]
Put together,
letting $R_4\ge R_3$ to be sufficiently large, we obtain that
$\tau_v \ge c_8 (\det(v))^{-1/2}$ with $c_8 = c_5/2$ for all $\|x\|
\ge
R_4$.

To sum up, by setting $C=\ball(0,R_4)$,
we get that for all $v\in\mathcal{C}^d$ with eigenvalues bounded from below
by $\kappa$, the estimate $P_v V(x) \le\lambda_v V(x)$ holds for
$x\notin
C$ with $\lambda_v \defeq1-c_8 \det(v)^{-1/2}$ satisfying
$(1-\lambda_v)^{-1} \le c_8^{-1}
\det(v)^{1/2}$. For $x\in{C}$, we have by (\ref{eq:drift-fraction}) that
$P_v V(x) \le2 V(x) \le2 \sup_{z\in C} V(z) \le b < \infty$, so
(\ref{eq:drift-ineq-m}) holds.
In the one-dimensional case, the above estimates can be applied separately
for the tails of the distribution.

Finally, set $\nu({B})\defeq|C|^{-1}|{B}\cap{C}|$,
and consider the minorization condition
(\ref{eq:mino-ineq-m}) for $x\in C$,
\begin{eqnarray*}
P_v(x,{B})
&\ge&\int_{B\cap C}
\biggl(1\wedge\frac{\pi(y)}{\pi(x)} \biggr) q_v(y-x) \,\ud y \\
&\ge&\frac{c_d}{\sqrt{\det(v)}}
\int_{B\cap C}
\biggl(1\wedge\frac{\pi(y)}{\pi(x)} \biggr)
\inf_{x,y\in C} e^{-{1}/{2}(x-y)v^{-1}(x-y)} \,\ud y\\
&\ge&\frac{c_d}{\sqrt{\det(v)}}
e^{-{1}/({2\kappa'})\operatorname{diam}(C)^2}
\frac{\inf_{z\in{C}}\pi(z)}{\sup_z\pi(z)}
\int_{B\cap C}
\,\ud y.
\end{eqnarray*}
So (\ref{eq:mino-ineq-m}) holds with $\delta_v \defeq c_9\det(v)^{-1/2}$
for some $c_9>0$. Finally, the claim holds with $c\defeq c_8^{-1}
\vee c_9^{-1}$.
\end{pf}
%

Finally, we are ready to prove the strong law of large numbers for the AM
process.
%
\begin{pf*}{Proof of Theorem \ref{thm:am-ergodicity}} 
We start by verifying the strong law of large numbers (\ref{eq:am-slln}).
Fix $t\ge1$ and consider first the constrained process
$(X_n^{(t)},S_n^{(t)})_{n\ge0}$ which is defined as the AM chain,
but with the constraint sets $K_n^{(t)}$ defined as
$K_n^{(t)} \defeq\{ s\in\mathbb{S}\dvtx|s|\le t n^{\varepsilon'} \}
$, with
$\varepsilon'=\varepsilon/(2d)$, and $\varepsilon\in(0,\kappa
_*^{-1}[(1/2)\wedge
(1-\alpha)])$, where $\kappa_*$ is the independent constant of
Theorem \ref{thm:slln-restricted}.

We check that assumptions \hyperlink{a:invariance}{(A1)}--\hyperlink{a:adapt-bound}{(A4)} are
satisfied
by the constrained process $(X_n^{(t)},S_n^{(t)})_{n\ge0}$ for all
$t\ge
1$. Condition
\hyperlink{a:invariance}{(A1)} is satisfied by construction of the Metropolis kernels
$P_s$.
Since $\det(v)\le\|v\|^d$, Proposition \ref{prop:geom-bound} ensures that
there is a compact $C\subset\R^d$ such that \hyperlink{a:uniform-drift-mino}{(A2)}
holds. For \hyperlink{a:kernel-lip}{(A3)}, we refer to \cite{andrieu-moulines},
Lemma 13, stating that $\|P_s f -P_{s'} f\|_{V^r} \le2
d\kappa^{-1} \|f\|_{V^r} |s^{(v)}-s'^{(v)}|$ for all
$s^{(v)},s'^{(v)}\in\mathcal{C}^d$ with eigenvalues bounded from
below by
$\kappa$.

Finally, we check that \hyperlink{a:adapt-bound}{(A4)} holds for any $\beta\in(0,1/2]$.
Similarly to \cite{sa-verifiable}, we have that
\begin{eqnarray*}
&&
\sup_{s\in K_n^{(t)}} \|H(s,x)\|_{V^\beta}\\
&&\qquad= \sup_{s\in K_n^{(t)}} \sup_{x\in\R^d}
\frac{|H(s,x)|}{V^\beta(x)} \\
&&\qquad\le\|\kappa I\| +
\sup_{x\in\R^d} \sup_{s\in K_n^{(t)}}
\frac{\|x\| + \|s^{(m)}\| + \|s^{(v)}\|
+ \| (x-s^{(m)})(x-s^{(m)})^T \|
}{V^\beta(x)} \\
&&\qquad\le\sqrt{d}\kappa+
\sup_{x\in\R^d}
\frac{\|x\| + \|x\|^2 + t^2n^{2\varepsilon'} + 2 t n^{\varepsilon'}
+ 2\|x\|t n^{\varepsilon'}
}{V^\beta(x)} \\
&&\qquad\le\sqrt{d}\kappa+
7t^2n^{2\varepsilon'}
\sup_{x\in\R^d}
\frac{\|x\|^2\vee1}{V^\beta(x)} \le\tilde{c} n^{\varepsilon}
\end{eqnarray*}
for any $\beta\in(0,1/2]$ by Lemma \ref{lemma:v-polybound}, where
$\tilde{c}
= \tilde{c}(t,\beta)$.
So, assumption \hyperlink{a:adapt-bound}{(A4)} holds for any $\beta\in
(0,1-\alpha)$.
In particular, we can select $\beta$ so that $\varepsilon< \kappa
_*^{-1}[(1/2)\wedge
(1-\alpha-\beta)]$. Clearly, $\sum_k k^{\kappa_*\varepsilon-1}
\eta_k
< \sum_k k^{\kappa_*\varepsilon- 2}<\infty$, so all the
conditions of Theorem \ref{thm:slln-restricted} are satisfied, implying
that the strong law of large numbers holds for the constrained process
$(X_n^{(t)},S_n^{(t)})$ for all $t\ge1$.

Define $B^{(t)}\defeq\{\forall n\ge0\dvtx S_n\in K_n^{(t)}\}$. We can
construct the constrained processes so that they coincide with the original
process in $B^{(t)}$. That is, for $\omega\in B^{(t)}$ we have
$(X_n(\omega),S_n(\omega)) = (X_n^{(t)}(\omega),S_n^{(t)}(\omega))$
for all
$n\ge0$. Lemma \ref{lemma:am-bound} ensures that we have $\P(\forall
n\ge
0\dvtx S_n\in K_n^{(t)})\ge g(t)$ where $g(t)\to1$ as $t\to\infty$.
As in the
proof of Theorem \ref{thm:slln-restricted}, we can use the Borel--Cantelli
lemma to deduce that (\ref{eq:am-slln}) holds almost surely.
\end{pf*}
%
%
\begin{remark} 
Since $\varepsilon>0$ can be selected arbitrarily small in the proof of
Theorem \ref{thm:am-ergodicity}, it is only required for (\ref{eq:am-slln})
to hold that the adaptation weights $\eta_n\in(0,1)$ are decreasing
and that
$\sum_k k^{\tilde{\varepsilon}-1} \eta_k <\infty$ holds for some
$\tilde{\varepsilon}>0$. In particular, one can choose $\eta_n
\defeq
(n+1)^{-\gamma}$ for any $\gamma>0$.
\end{remark}
%
%
\begin{remark} 
Condition (\ref{eq:super-pexp}) implies the super-exponential
decay of the tails of $\pi$:
%
%
\begin{equation} \label{eq:super-exp}
\lim_{r\to\infty} \sup_{\|x\|\ge r}
\frac{x}{\|x\|} \cdot\nabla\log\pi(x) = -\infty.
\end{equation}
This condition, with the contour regularity condition (\ref{eq:reg-contour}),
are common conditions to
ensure geometric ergodicity of a random-walk Metropolis algorithm, and
many standard distributions fulfil them
\cite{jarner-hansen}. The decay
condition (\ref{eq:super-pexp}) is only slightly more stringent than
(\ref{eq:super-exp}).
\end{remark}
%

Finally, we formulate a central limit theorem for the AM algorithm.
%
\begin{theorem}
\label{thm:am-clt} 
Assume $\pi$ satisfies the conditions of Theorem \ref{thm:am-ergodicity}.
For any $f$ with $\|f\|_{V^\alpha} < \infty$ for some
$0\le\alpha< 1/2$, where
$V(x)\defeq c_V \pi^{-1/2}(x)$ and $c_V = (\sup_x \pi(x))^{1/2}$,
it holds that
\[
\frac{1}{\sqrt{n}}\sum_{k=1}^n
[f(X_k)-\pi(f) ] \stackrel{n\to\infty}{\hbox to 1cm{\rightarrowfill}}
N(0,\sigma^2)
\]
in distribution, where $\sigma^2\in[0,\infty)$ is a constant.
\end{theorem}

The proof of Theorem \ref{thm:am-clt} follows
by the techniques of the present paper applied to \cite
{andrieu-moulines}, Theorem 9.
A fully detailed proof can be found in the preprint
\cite{saksman-vihola-preprint}.


\begin{appendix}
\section{\texorpdfstring{Proof of Lemma
\protect\lowercase{\ref{lemma:drift-to-conv}}}{Proof of Lemma 3}}
\label{sec:proof-lemma-drift-to-conv} 

We provide a restatement of a part of a theorem by
Meyn and Tweedie \cite{meyn-tweedie-computable}
before proving Lemma
\ref{lemma:drift-to-conv}. For a more recent work on quantitative
convergence bounds, we refer to \cite{baxendale-bounds}.
%
\begin{theorem}
\label{th:m-t-computable} 
Suppose that the following drift and minorization conditions hold:
\begin{eqnarray*}
P V(x) &\le& \lambda V(x) + b \mathbh{1}_{C}(x)\qquad
\forall x\in\mathbb{X}, \\
P(x,A) &\ge& \delta\nu(A)\qquad \forall x\in C, \forall
{A}\subset\mathbb{X}
\end{eqnarray*}
for constants $\lambda<1$, $b<\infty$ and $\delta>0$,
a set ${C}\subset\mathbb{X}$ and a
probability measure $\nu$ on ${C}$.
Moreover, suppose that $\sup_{x\in{C}}V(x)\le b$.
Then, for all $k\ge1$,
\[
\|P_s^k(x,\uarg)-\pi(\uarg)\|_V \le V(x) (1+\gamma)
\frac{\rho}{\rho-\vartheta}\rho^k
\]
for any $\rho>\vartheta=1-\tilde{M}^{-1}$, for
\[
\tilde{M} =
\frac{1}{(1-\check{\lambda})^2} \bigl[1-\check{\lambda}+\check
{b}+\check{b}^2
+ \bar{\zeta} \bigl( \check{b}(1-\check{\lambda}) \check{b}^2 \bigr) \bigr]
\]
defined in terms of
\begin{eqnarray*}
\gamma&=& \delta^{-2} [4b+2\delta\lambda b ] ,\\
\check{\lambda} &=& (\lambda+\gamma)/(1+\gamma)<1, \\
\check{b} &=& b + \gamma< \infty
\end{eqnarray*}
and the bound
\[
\bar{\zeta} \le
\frac{4-\delta^2}{\delta^5} \biggl(\frac{b}{1-\lambda} \biggr)^2.
\]
\end{theorem}
%
%
\begin{pf}
See \cite{meyn-tweedie-computable}, Theorem 2.3.
\end{pf}
%
%
\begin{pf*}{Proof of Lemma \ref{lemma:drift-to-conv}} 
Observe that $P_s V(x) = \mathbb{E}[V(X_{n+1}) | X_n=x,S_n=s]$, and therefore
by Jensen's inequality, \hyperlink{a:uniform-drift-mino}{(A2)} implies
for $x\notin{C}_n$ that
\[
P_s V^r(x) \le(P_s V(x))^r \le\lambda_n^r V^r(x).
\]
We can bound $\tilde{\lambda}_n\defeq\lambda_n^r \le(1-c^{-1}
n^{-\varepsilon})^r
\le1 - r c^{-1}n^{-\varepsilon}$ implying
\[
(1-\tilde{\lambda}_n)^{-1} \le r^{-1} c n^{\varepsilon},
\]
whenever $r\in(0,1]$. Similarly, for $x\in{C}_n$, one has
$P_s V^r(x) \le(\sup_{z\in{C}_n}V(z) + b_n)^r \le(2b_n)^r$, so by
letting $\tilde{b}_n \defeq(2b_n)^r$,
we obtain the drift inequality
\[
P_s V^r(x) \le\tilde{\lambda}_n V^r(x) + \tilde{b}_n\mathbh{1}_{C_n}(x),
\]
and we can bound $\tilde{b}_n \le(2c n^\varepsilon)^r$. We have the bound
$(1-\tilde{\lambda}_n)^{-1}\vee\tilde{b}_n \le\tilde{c}
n^{\varepsilon}$
with some $\tilde{c}=\tilde{c}(c,r)\ge1$.

Now, we can apply Theorem \ref{th:m-t-computable}, where we can
estimate the
constants
\begin{eqnarray*}
\gamma_n &=& \delta_n^{-2} [4\tilde{b}_n+2\delta_n\tilde{\lambda}_n
\tilde{b}_n ]
\le(c n^{\varepsilon})^2 6(\tilde{c}n^\varepsilon) = a_1
n^{3\varepsilon}, \\
\check{b}_n &=& \tilde{b}_n + \gamma_n
\le(\tilde{c}+a_1)n^{3\varepsilon}
\le a_2 n^{3\varepsilon}
\end{eqnarray*}
and consequently
\[
1-\check{\lambda}_n
= \frac{1-\tilde{\lambda}_n}{1+\gamma_n}
\ge\frac{\tilde{c}^{-1} n^{-\varepsilon}}{1+a_1 n^{3\varepsilon}}
\ge\frac{\tilde{c}^{-1}}{1+a_1} n^{-4\varepsilon} = a_3^{-1}
n^{-4\varepsilon}.
\]
Moreover,
\[
\bar{\zeta}_n \le
\frac{4-\delta_n^2}{\delta_n^5} \biggl(\frac{\tilde{b}_n}{1-\tilde
{\lambda
}_n} \biggr)^2
\le4 (c n^{\varepsilon})^5(\tilde{c}n^\varepsilon)^2 (\tilde
{c}n^{\varepsilon})^2
= a_4 n^{9\varepsilon},
\]
and then
\begin{eqnarray*}
\tilde{M}_n &=& \frac{1}{(1-\check{\lambda})^2}
\bigl[ 1 - \check{\lambda}_n + \check{b}_n + \check{b}_n^2
+ \bar{\zeta}_n\bigl(\check{b}_n(1-\check{\lambda}_n) + \check{b}_n^2\bigr)
\bigr]\\
&\le& (a_3 n^{4\varepsilon})^2
[ 1 + \check{b}_n + \check{b}_n^2
+ \bar{\zeta}_n(\check{b}_n + \check{b}_n^2) ] \\
&\le& (a_3 n^{4\varepsilon})^2(5 \bar{\zeta}_n\check{b}_n^2)
\le5 a_3^2 n^{8\varepsilon} a_4 n^{9\varepsilon} a_2^2
n^{6\varepsilon
} = a_5
n^{23\varepsilon}
\end{eqnarray*}
since we can assume that $\check{b}_n,\bar{\zeta}_n\ge1$. Now,
\[
1-\vartheta_n = \tilde{M}_n^{-1} \ge a_5^{-1} n^{-23\varepsilon}
\]
and we can choose $\rho_n \in(\vartheta_n,1)$ by letting
$\rho_n \defeq\frac{1+\vartheta_n}{2}$. We have
\[
\rho_n - \vartheta_n
= 1-\rho_n
= \tfrac{1}{2}(1-\vartheta_n)
\ge\tfrac{1}{2}c_{9}^{-1} n^{-23\varepsilon} = (a_6 n^{23\varepsilon})^{-1}.
\]
Finally, from Theorem \ref{th:m-t-computable}, one obtains the bound
\[
{\|P_s^k(x,\uarg)-\pi(\uarg)\|}_{V^r} \le V^r(x) L_n \rho_n^k,
\]
where
\begin{eqnarray*}
(1-\rho_n)^{-1} &\le& a_6 n^{23\varepsilon}, \\
L_n &=& (1+\gamma_n)\frac{\rho_n}{\rho_n-\vartheta_n}
\le(1+a_1 n^{3\varepsilon})(a_6 n^{23\varepsilon})
\le a_7 n^{26\varepsilon}
\end{eqnarray*}
with $a_7 = (1+a_1)a_6$. This concludes the proof with
$\kappa_2 = 26$ and $c_2 = a_7$.
\end{pf*}
%

\vspace*{-12pt}

\section{Birnbaum and Marshall's inequality}
\label{sec:martingale} 

%
\begin{theorem}[(Birnbaum and Marshall)]
\label{th:birnbaum-marshall} 
Let $(X_k)_{k=1}^n$ be random variables, such that
\[
\mathbb{E}[|X_k| | \F_{k-1}]
\ge\psi_k|X_{k-1}|,
\]
where $\F_k\defeq\sigma(X_1,\ldots,X_k)$, and $\psi_k\ge0$.
Let $a_k > 0$, and define
\[
b_k \defeq\max\Biggl\{
a_k, a_{k+1}\psi_{k+1}, \ldots, a_n \prod_{j=k+1}^n \psi_j
\Biggr\}
\]
for $1\le k\le n$, and $b_{n+1}\defeq0$.
If $p\ge1$ is such that $\E|X_k|^p<\infty$ for all $1\le k \le n$, then
\[
\P\Bigl(\max_{1\le k\le n} a_k|X_k| \ge1 \Bigr)
\le\sum_{k=1}^n (b_k^p - \psi_{k+1}^p b_{k+1}^p)\E|X_k|^p.
\]
\end{theorem}
%
%
\begin{pf} 
See \cite{birnbaum-marshall}, Theorem 2.1.
\end{pf}
%
%
%
\begin{corollary} 
\label{cor:birnbaum-marshall-martingale}
Let $(M_k)_{k=1}^n$ be a martingale with respect to
$(\F_k)_{k=1}^n$. Let $(a_k)_{k=1}^n$ be a
strictly positive nonincreasing sequence. If $p\ge1$ is such that
$\E|M_k|^p<\infty$ for all $1\le k\le n$, then for $1\le m\le n$,
\[
\P\Bigl( \max_{m\le k\le n} a_k|M_k| \ge1 \Bigr)
\le a_{n}^p \E|M_n|^p + \sum_{k=m}^{n-1} (a_k^p - a_{k+1}^p)\E|M_k|^p.
\]
\end{corollary}
%
%
\begin{pf} 
By Jensen's inequality,
\[
\mathbb{E}[|M_k| | \F_{k-1}] \ge| \mathbb{E}[M_k | \F_{k-1}] |
= |M_{k-1}|.
\]
Define $\psi_k \defeq1$
for $1\le k\le n$, and $\tilde{a}_k \defeq a_m$
for $1\le k \le m$ and $\tilde{a}_k \defeq a_k$ for $m< k\le n$.
The result follows from Theorem \ref{th:birnbaum-marshall}.\vadjust{\goodbreak}
\end{pf}
%


\section{Contour surface containment}
\label{sec:contour} 

%
\begin{lemma}
\label{lemma:contour} 
Suppose $A\subset\R^d$ is a smooth surface parameterized by the unit sphere
$\mathcal{S}^d$, that is,
$A = \{ug(u)\dvtx u\in\mathcal{S}^d\}$ with a continuously differentiable
radial function
$g\dvtx\mathcal{S}^d\to(0,\infty)$.
Assume also that outer-pointing normal $n$ of $A$ satisfies
$n(x) \cdot x/\|x\| \ge\beta$ for all $x\in A$ with some constant
$\beta> 0$.
There is a constant $M<\infty$ depending only on $\beta$ such that for
any $x,y\in A$, it holds that
$M^{-1} \le\|x\|/\|y\| \le M$.
\end{lemma}
\begin{pf} 
Consider first the two-dimensional case.
Let $x$ and $y$ be two distinct points in $A$.
We employ polar coordinates, thus let $u(\theta)r(\theta)\in A$
with $u(\theta) \defeq[\cos(\theta),\sin(\theta)]^T$ and
$r(\theta)
\defeq
g(u(\theta))$
so that $u(\theta_1)r(\theta_1)=x$ and $u(\theta_2)r(\theta_2)=y$
with $\theta_1,\theta_2\in[0,2\pi)$.

Let $\alpha(\theta)$ stand for the (smaller) angle between
$u(\theta)$ and the normal of the curve~$A$, that is, the curve parametrized
by $\theta\to u(\theta)r(\theta)$. Our assumption says that
$|\alpha(t)|\leq\alpha_0 \defeq\arccos(\beta)<\pi/2$ for all
$\theta\in[0,2\pi]$.
On the other hand, an elementary computation shows that
\[
\tan(\alpha(\theta) )=\frac{r'(\theta)}{r(\theta)},
\]
and hence we have $|\frac{\ud}{\ud\theta}\log
r(\theta))|=|r'(\theta)/r(\theta)|\leq\tan\alpha_0$ uniformly.
We may estimate $|{\log}\|x\| - {\log}\|y\|| \le2\pi\tan(\alpha_0)$
yielding the claim with $M = e^{2\pi\tan{\alpha_0}}$.

For $d\ge3$,
take the plane $T$ containing the origin and the points
$x$ and $y$. This reduces the situation to two dimensions,
since $A\cap T$ inherits the given normal condition
of the surface and the radius vector.
\end{pf}
\end{appendix}


\section*{Acknowledgments} 

We thank the anonymous referees for a careful review and comments improving
the paper significantly. We also thank Gersende Fort for useful comments.



%
\printaddresses

\end{document}